\documentclass{article}

\usepackage{hyperref}

\usepackage{hyperref}

\usepackage{amssymb}

\newcommand{\newc}{\newcommand}

\newc{\eqnoset}{\setcounter{equation}{0}}

\newcommand{\mref}[1]{(\ref{#1})}
\newcommand{\reflemm}[1]{Lemma~\ref{#1}}
\newcommand{\refrem}[1]{Remark~\ref{#1}}
\newcommand{\reftheo}[1]{Theorem~\ref{#1}}

\newcommand{\refcoro}[1]{Corollary~\ref{#1}}

\newcommand{\refsec}[1]{Section~\ref{#1}}

\newcommand{\beq}{\begin{equation}}
	\newcommand{\eeq}{\end{equation}}
\newcommand{\beqno}[1]{\begin{equation}\label{#1}}
	
	\newcommand{\barr}{\begin{array}}
		\newcommand{\earr}{\end{array}}
	
	\newc{\bearr}{\begin{eqnarray*}}
		\newc{\eearr}{\end{eqnarray*}}
	
	\newc{\bearrno}[1]{\begin{eqnarray}\label{#1}}
		\newc{\eearrno}{\end{eqnarray}}
	
	\newc{\non}{\nonumber}
	\newc{\nol}{\nonumber\nl}
	
	\newcommand{\bdes}{\begin{description}}
		\newcommand{\edes}{\end{description}}
	\newc{\benu}{\begin{enumerate}}
		\newc{\eenu}{\end{enumerate}}
	\newc{\btab}{\begin{tabular}}
		\newc{\etab}{\end{tabular}}

	\newtheorem{theorem}{Theorem}[section]
	\newtheorem{defi}[theorem]{Definition}
	\newtheorem{lemma}[theorem]{Lemma}
	\newtheorem{rem}[theorem]{Remark}
	\newtheorem{exam}[theorem]{Example}
	\newtheorem{propo}[theorem]{Proposition}
	\newtheorem{corol}[theorem]{Corollary}
	\newtheorem{conj}[theorem]{Conjecture}

	\newcommand{\btheo}[1]{\begin{theorem}\label{#1}}
		\newc{\brem}[1]{\begin{rem}\label{#1}\em}
			\newc{\bexam}[1]{\begin{exam}\label{#1}\em}
				\newc{\bdefi}[1]{\begin{defi}\label{#1}}
					\newcommand{\blemm}[1]{\begin{lemma}\label{#1}}
						\newcommand{\bprop}[1]{\begin{propo}\label{#1}}
							\newcommand{\bcoro}[1]{\begin{corol}\label{#1}}
								\newcommand{\btheoc}[1]{\begin{conj}\label{#1}}
									\newcommand{\etheo}{\end{theorem}}
								\newc{\etheoc}{\end{conj}}
							\newcommand{\elemm}{\end{lemma}}
						\newcommand{\eprop}{\end{propo}}
					\newcommand{\ecoro}{\end{corol}}
				\newc{\erem}{\end{rem}}
			\newc{\eexam}{\end{exam}}
		\newc{\edefi}{\end{defi}}
	
	\newc{\rmk}[1]{{\bf REMARK #1: }}
	\newc{\DN}[1]{{\bf DEFINITION #1: }}
	
	\newcommand{\bproof}{{\bf Proof:~~}}
	\newc{\eproof}{{\vrule height8pt width5pt depth0pt}\vspace{3mm}}
	
	\newc{\bfrac}[2]{\dspl{\frac{#1}{#2}}}

	\newc{\nid}{\noindent}
	
	\newcommand{\dspl}{\displaystyle}
	\newc{\grad}{\nabla}
	\newc{\Div}{\mbox{div}}
	\newc{\pdt}[1]{\dspl{\frac{\partial{#1}}{\partial t}}}
	\newc{\pdn}[1]{\dspl{\frac{\partial{#1}}{\partial \nu}}}
	\newc{\pdNi}[1]{\dspl{\frac{\partial{#1}}{\partial \mathcal{N}_i}}}
	\newc{\pD}[2]{\dspl{\frac{\partial{#1}}{\partial #2}}}
	\newc{\dt}{\dspl{\frac{d}{dt}}}
	\newc{\bdry}[1]{\mbox{$\partial #1$}}
	\newc{\sgn}{\mbox{sign}}
	
	\newc{\Hess}[1]{\frac{\partial^2 #1}{\pdh z_i \pdh z_j}}
	\newc{\hess}[1]{\partial^2 #1/\pdh z_i \pdh z_j}

	\newc{\ag}{\alpha}
	\newc{\bg}{\beta}
	\newc{\cg}{\gamma}\newc{\Cg}{\Gamma}
	\newc{\dg}{\delta}\newc{\Dg}{\Delta}
	\newc{\eg}{\varepsilon}
	\newc{\zg}{\zeta}
	\newc{\thg}{\theta}
	\newc{\llg}{\lambda}\newc{\LLg}{\Lambda}
	\newc{\kg}{\kappa}
	\newc{\rg}{\rho}
	\newc{\sg}{\sigma}\newc{\Sg}{\Sigma}
	\newc{\tg}{\tau}
	\newc{\fg}{\phi}\newc{\Fg}{\Phi}
	\newc{\vfg}{\varphi}
	\newc{\og}{\omega}\newc{\Og}{\Omega}
	\newc{\pdh}{\partial}
	
	\newc{\ccG}{{\cal G}}

	\newc{\ii}[1]{\int_{#1}}
	\newc{\iidx}[2]{{\dspl\int_{#1}~#2~dx}}
	\newc{\bii}[1]{{\dspl \ii{#1} }}
	\newc{\biii}[2]{{\dspl \iii{#1}{#2} }}
	\newc{\su}[2]{\sum_{#1}^{#2}}
	\newc{\bsu}[2]{{\dspl \su{#1}{#2} }}

	\newc{\biiom}[1]{{\dspl\int_{\bdrom}~ #1 ~d\sg}}
	\newc{\io}[1]{{\dspl\int_{\Og}~ #1 ~dx}}
	\newc{\bio}[1]{{\dspl\int_{\bdrom}~ #1 ~d\sg}}
	\newc{\bsir}{\bsu{i=1}{r}}
	\newc{\bsim}{\bsu{i=1}{m}}
	
	\newc{\iibr}[2]{\iidx{\bprw{#1}}{#2}}
	\newc{\Intbr}[1]{\iibr{R}{#1}}
	\newc{\intbr}[1]{\iibr{\rg}{#1}}
	\newc{\intt}[3]{\int_{#1}^{#2}\int_\Og~#3~dxdt}
	
	\newc{\itQ}[2]{\dspl{\int\hspace{-2.5mm}\int_{#1}~#2~dz}}
	\newc{\mitQ}[2]{\dspl{\rule[1mm]{4mm}{.3mm}\hspace{-5.3mm}\int\hspace{-2.5mm}\int_{#1}~#2~dz}}
	\newc{\mitQQ}[3]{\dspl{\rule[1mm]{4mm}{.3mm}\hspace{-5.3mm}\int\hspace{-2.5mm}\int_{#1}~#2~#3}}
	
	\newc{\mitx}[2]{\dspl{\rule[1mm]{3mm}{.3mm}\hspace{-4mm}\int_{#1}~#2~dx}}
	\newc{\mitmu}[2]{\dspl{\rule[1mm]{3mm}{.3mm}\hspace{-4mm}\int_{#1}~#2~d\mu}}
	\newc{\iidmu}[2]{\iidx{#1}{#2}}
	
	\newc{\iidm}[3]{{\dspl\int_{#1}~#2~d #3}}
	
	\newc{\itQmu}[2]{\dspl{\int\hspace{-2.5mm}\int_{#1}~#2~d\mu}}
	\newc{\mitQmu}[2]{\dspl{\rule[1mm]{4mm}{.3mm}\hspace{-5.3mm}\int\hspace{-2.5mm}\int_{#1}~#2~d\mu}}

	\newc{\mitQq}[2]{\dspl{\rule[1mm]{4mm}{.3mm}\hspace{-5.3mm}\int\hspace{-2.5mm}\int_{#1}~#2~d\bar{z}}}
	\newc{\itQq}[2]{\dspl{\int\hspace{-2.5mm}\int_{#1}~#2~d\bar{z}}}

	\newc{\pder}[2]{\dspl{\frac{\partial #1}{\partial #2}}}

	\newc{\bdrom}{\bdry{\Og}}

	\newc{\bilhom}{\mbox{Bil}(\mbox{Hom}(\RR^{nm},\RR^{nm}))}
	\newc{\VV}[1]{{V(Q_{#1})}}
	
	\newc{\ccA}{{\mathcal A}}
	\newc{\ccB}{{\mathcal B}}
	\newc{\ccC}{{\mathcal C}}
	\newc{\ccD}{{\mathcal D}}
	\newc{\ccE}{{\mathcal E}}
	\newc{\ccH}{\mathcal{H}}
	\newc{\ccF}{\mathcal{F}}
	\newc{\ccI}{{\mathcal I}}
	\newc{\ccJ}{{\mathcal J}}
	\newc{\ccK}{{\mathcal K}}
	\newc{\ccP}{{\mathcal P}}
	\newc{\ccQ}{{\mathcal Q}}
	\newc{\ccR}{{\mathcal R}}
	\newc{\ccS}{{\mathcal S}}
	\newc{\ccT}{{\mathcal T}}
	\newc{\ccX}{{\mathcal X}}
	\newc{\ccY}{{\mathcal Y}}
	\newc{\ccZ}{{\mathcal Z}}
	
	\newc{\bb}[1]{{\mathbf #1}}
	
	\newc{\myprod}[1]{\langle #1 \rangle}
	\newc{\mypar}[1]{\left( #1 \right)}
	
	\newc{\BLLg}{\mathbf{\LLg}}
	
	\newc{\mA}{\mathbf{A}}
	\newc{\mB}{\mathbf{B}}
	\newc{\mC}{\mathbf{C}}
	\newc{\mD}{\mathbf{D}}
	\newc{\mE}{\mathbf{E}}
	\newc{\mF}{\mathbf{F}}
	\newc{\mJ}{\mathbf{J}}
	\newc{\mG}{\mathbf{G}}
	\newc{\mP}{\mathbf{P}}
	\newc{\mR}{\mathbf{R}}
	\newc{\mQ}{\mathbf{Q}}
	\newc{\mX}{\mathbf{X}}
	\newc{\muu}{\mathbf{u}}
	\newc{\mvv}{\mathbf{v}}
	
	\newc{\mllg}{\mathbb{\lambda}}
	\newc{\mLLg}{\mathbf{\LLg}}

	\newc{\lspn}[2]{\mbox{$\| #1\|_{\Lsp{#2}}$}}
	\newc{\Lpn}[2]{\mbox{$\| #1\|_{#2}$}}
	\newc{\Hn}[1]{\mbox{$\| #1\|_{H^1(\Og)}$}}

	\newc{\mynorm}[2]{\| #1\|_{#2}}

	\newcommand{\RR}{{\rm I\kern -1.6pt{\rm R}}}

	\newc{\itQQ}[2]{\dspl{\int_{#1}#2\,dz}}
	\newc{\mmitQQ}[2]{\dspl{\rule[1mm]{4mm}{.3mm}\hspace{-4.3mm}\int_{#1}~#2~dz}}
	\newc{\MmitQQ}[2]{\dspl{\rule[1mm]{4mm}{.3mm}\hspace{-4.3mm}\int_{#1}~#2~d\mu}}

	\newc{\MUmitQQ}[3]{\dspl{\rule[1mm]{4mm}{.3mm}\hspace{-4.3mm}\int_{#1}~#2~d#3}}
	\newc{\MUitQQ}[3]{\dspl{\int_{#1}~#2~d#3}}

	\newc{\mccP}{\mathbb{P}}
	\newc{\mccK}{\mathbb{K}}
	
	\newc{\DKTmU}{\mccK(U)}
	\newc{\DKTmUold}{(K_U(U)^{-1})^T}
	
	\newc{\myPi}{\mathbf{W}}
	\newc{\myIbar}{\bar{\ccI}_1}
	\newc{\myIhat}{\hat{\ccI}_1}
	\newc{\myIbreve}{\breve{\ccI}_0}
	
	\newc{\mmk}{\mathbf{k}}

	\newcommand{\ma}{\mathbf{a}}
	\newcommand{\mg}{\mathbf{g}}

	\newc{\mfu}{\mathbf{f_u}}
	\newc{\mh}{\mathbf{h}}
	\newc{\mb}{\mathbf{b}}
	\newc{\mf}{\mathbf{f}}

	\newcommand{\barrl}[2]{\barr{ll}\lefteqn{#1}\hspace{#2}&\\}

	\newc{\twomatrix}[1]{\left[\barr{cc}#1\earr\right]}
	\newc{\threematrix}[1]{\left[\barr{ccc}#1\earr\right]}

	\newc{\mN}{\mathbf{N}}
	\newc{\mI}{\mathbf{I}}
	\newc{\mH}{\mathbf{H}}

	\newc{\mk}{\mathbf{k}}
	\newc{\mr}{\mathbf{r}}

	\newc{\DIAGM}[2]{\left[\barr{ccc}#1&0\ldots&0\\
		\vdots&\ddots&\vdots\\0&\ldots0&#2\earr \right]}
	\newc{\DiagM}[2]{\mbox{diag}\left[#1
		\cdots #2 \right]}
	\newc{\vVEC}[2]{\left[\barr{c}#1\\
		\vdots\\#2\earr \right]}
	\newc{\hVEC}[2]{\left[#1
		\cdots #2 \right]}
	
	\newc{\mq}{\mathbf{q}}

	\newc{\msys}[1]{\left\{\barr{l}#1\earr
		\right.}
	\newc{\msysa}[1]{\left\{\barr{ll}#1\earr
		\right.}
	
	\newc{\bbM}{\mathbb{M}}
	\newc{\mat}[1]{\left[\barr{cc}#1\earr\right]}
	
	\newc{\me}{\mathbf{e}}

	\newc{\vecc}[2]{\left[\barr{cc}#1\\#2\earr\right]}
	\newc{\mL}{\mathbb{L}}
	
	\newc{\cO}{{\cal O}}
	
	\newc{\cM}{{\cal M}}
	
	\newc{\myega }{\eg_0(R)}
	\newc{\myeg}{\eg_1(\eg_*)}
	\newc{\myegp}{\hat{\eg}_1(\eg_*)}

\newc{\diagA}{\mathbb{A}_d}
\newc{\mBB}{\mathbb{B}}
\newc{\MLT}[1]{{\cal M}_{lt}(\Og,#1)}
\newc{\ALT}[1]{{\cal A}_{l}(\Og,#1)}
\newc{\mM}{\mathbb{M}}

\newc{\diag}[1]{\mbox{diag}(#1)}
\newc{\off}[1]{\mbox{offdiag}(#1)}
\newc{\mT}{\mathbb{T}}

\usepackage{amssymb}

\begin{document}

	\vspace*{-.8in}
	\begin{center} {\LARGE\em A priori Bound of Solutions to a Class of Elliptic/Parabolic Cross Diffusion Systems of $m$ Equations on Two/Three Dimensional Domains. Existence on thin domains}
		
	\end{center}

	\vspace{.1in}
	
	\begin{center}

		{\sc Dung Le}{\footnote {Department of Mathematics, University of
				Texas at San
				Antonio, One UTSA Circle, San Antonio, TX 78249. {\tt Email: Dung.Le@utsa.edu}
				}}

	\end{center}

	\begin{abstract}
		We establish several bounds for solutions to elliptic/parabolic cross-diffusion systems of $m$ equations ($m\ge2$) on 2d/3d domains $\Og$. We settle the existence and global existence problems in these cases and also provide new counterexamples for the case of general dimensions. Most importantly, we prove that when $m=N=3$,  the thinness of $\Og$ in $\RR^3$ is sufficient and necessary. When $m$ is arbitrary and $N=3$, we establish global existence results for nonlinear cross-diffusion systems (the case of scalar  semilinear equations was considered in the literature but the classical  methods do not seem to be applicable here).
		\end{abstract}

\section{Introduction}

This paper deals with two cross diffusion systems, the elliptic one
\beqno{exsysz}-\Div(A(W)DW)=G(W)W\eeq
and its corresponding evolution 
\beqno{exsyszpara}W_t=\Div(A(W)DW)+G(W)W\eeq
on bounded domain $\Og$ of $\RR^N$ with homogeneous Dirichlet or Neumann boundary conditions. Here, $W$ is an unknown vector in $\RR^m$, $m\ge 2$. As usual, $A(W), G(W)$ are $m\times m$ matrices and for some $C,k\ge0, c>0$ we assume the ellipticity and growth conditions
$$\myprod{A(W)DW,DW}\ge c(1+|W|^k)|DW|^2\mbox{ and } |G(W)|\le C(|W|^k+1).$$

We will present several a priori estimates for higher norms strong and weak solutions of the above two systems by weaker norms under these conditions. These estimates play crucial roles in the existence of problem \mref{exsysz} and global existence of problem \mref{exsyszpara}.

The existence of problem \mref{exsysz} can be solved by using the fixed point index theory (see e.g. \cite{Amsurvey} for diagonal systems and \cite{dlebook} for nondiagonal cases) and of course, one must establish that all fixed points of certain compact map associated to the system have uniform bound of the norm in the space where the map was defined. To have 
the compactness, the spaces have been usually H\"older spaces so that the approaches frequently faced regularity problems, a hard and opened fundamental one in the theory of partial differential systems like \mref{exsysz} when $A$ is full. 

We have seen (\cite{Gius} and recently \cite{dlejfa,dlebook}) the regularity theory for elliptic systems has a closed tie will the BMO norms of solutions in small balls, which is somewhat easier but still a very hard problem by no mean. However, when $N=2$, the problem can be done by estimating $\|DW\|_{L^2(\Og)}$. To gain some insights, in \refsec{uniellsec} we consider the easiest case of elliptic systems of 2 equations on planar domains ($m=N=2$). Combine with the dual methods and \cite{lehaljde}, we will estimate $\|DW\|_{L^2(\Og)}$ in terms of the weakest norm $\|W\|_{L^1(\Og)}$, which can be controlled by looking at the algebraic structure of $A,G$ of \mref{exsysz}. The uniform estimate for $C^1$ norms of solutions is a difficult problem can there is a vast literature on this. Here, we mention a minimum requirement. Assume that the solutions are nonnegative. 
We will extend the usual interpolation inequality in Sobolev spaces to new nonlinear ones which involved with the principal eigenpairs of $A$, which exist via the theory of cooperative systems (e.g. \cite{lopez}). The key point is that when {\em $|W|$ large and $G(W)$ is cooperative but we can bound the solutions uniformly}. New results on coexistence (existence of nontrivial solutions) will br reported in a forthcoming paper (\cite{lecoex}).

The above use of principal eigenpairs does not exist (yet) for parabolic case (unless they are periodic). In \refsec{uniparasec}, we consider the global existence and dynamics of solutions to \mref{exsyszpara} and deal with their  estimates and regularity. The problems are even harder and again, thanks to the pioneering works (e.g. \cite{GiaS,Am2}) we see that one has to  show that BMO norms of solutions are small in small balls. The later was opened and there are counterexamples when $m=N\ge 3$ (see \cite{JS} and reference therein). Here, we allow $m\ge2$ is arbitrary and settle the problem when $N=2$ and also provide new blow-up examples (for nonlinear systems for general $m,N\ge2$) to see that our hypotheses are almost optimal.

Importantly, we report global existence results when $N=3$ ($m$ is arbitrary) on thin domains. These are new when $A$ is a full matrix and nonlinear (the case \mref{exsyszpara} is scalar and semilinear was considered in the literature but their methods do not seem to be applicable here). Combine with the counterexamples in \cite{JS}, we see that the thinness of $\Og$ in $\RR^3$ is sufficient and necessary.

	\section{The elliptic case:}\label{uniellsec} \eqnoset
	For simplicity we consider the case of 2 equations for $W=(u,v)$. Let $A(W), G(W)$ be $2\times 2$ matrices and for some $C,k\ge0, c>0$ we have \beqno{gencond}\myprod{A(W)DW,DW}\ge c(1+|W|^k)|DW|^2\mbox{ and } |G(W)|\le C(|W|^k+1). \eeq Testing the system $-\Div(A(W)DW)=G(W)W$ with $W$, we obtain
	\beqno{keyboundest}c\iidx{B}{(1+|W|^k)|DW|^2}\le C\iidx{B}{(1+|W|^k)|W|^2}.\eeq
	
	If we can absorb the integral on the right into that of $|W|^k|DW|^2$ on the left then we obtain a uniform bound for $\|DW\|_{L^2(B)}$
	so that the $BMO$ of weak solutions are uniformly small in small balls if $N=2$. We then have that $C^{\ag}$ norm of solutions for some $\ag>0$ are uniformly bounded by the theory in \cite{dlebook1}. This provides the existence of a solution of the system.

	Consider the simplest case. Let $A$ be a constant matrix and $G(W)$ be bounded (so $k=0$). It is well known that for any $\eg>0$ there is $C(\eg)$ such that \beqno{intineq0}\iidx{B}{|W|^2}, \|W\|_{L^{2^*}(B)}^2\le \eg\iidx{B}{|DW|^2}+C(\eg)\left(\iidx{B}{W}\right)^2,\eeq
	so the uniform bound for $\|W\|_{L^1(B)}$  yields that  of $\|DW\|_{L^2(B)}$ (see \cite[Lemma 5.3.5]{dlebook1}).

	Denote by $(\llg,\fg)$ the principal eigenpairs of $-\Delta$. Testing \mref{exsysz} with $\Fg=\mat{\fg\\\fg}$ and writing $x:=\myprod{W,\Fg}_{L^2}$, we get $\llg x=A^{-1}\myprod{G(W)W,\Fg}_{L^2}$. Since $W,\fg>0$ in $B$, $x\in\RR^m$ is positive and, obviously, $\|W\|_{L^1}$ is bounded if $x$ is. Since $x>0$ the latter is equivalent to the fact that there are $x_0,n_0$ in the positive cone of $\RR^m$ such that $\myprod{\llg x-x_0,n_0}\le0$ or 
	$$\myprod{\myprod{A^{-1}G(W)W,\Fg}_{L^2}-x_0,n_0}\le0\Leftrightarrow\myprod{\myprod{G(W)W,\Fg}_{L^2}-x_1,z_0}\le0$$ where $x_1=A x_0$ and $z_0=(A^{-1})^Tn_0$ for some $x_0,n_0>0$. As $\fg>0$ and we can normalize it such that $\|\fg\|_{L^\infty}(B)=1$, the above holds if $$\iidx{B}{\myprod{z_0,G(W)W\Fg}}\le C_0\Leftrightarrow\myprod{z_0,G(w)w}\le c_0\quad\mbox{for some $c_0\in L^1(B)$ and all $w\ge0$}.$$ Note that $z_0=(A^{-1})^Tn_0$ for some $n_0>0$ but $z_0$ does not have to be a positive vector.

	The above argument can be extended to the nonconstant $A$ given by
	$$A(u,v)=[\ma^{ij}]=\mat{a_{11}\cg_1(u)&a_{12}\cg_2(v)\\a_{21}\cg_1(u)&a_{22}\cg_2(v)},$$ where $a_{ij}$'s are constants and $\cg_i$'s are positive functions. Assume that $[a_{ij}]$ is elliptic. It is easy to see that the above calculation can be repeated to give that $\|\Cg(W)\|_{L^1(B)}$'s are uniformly bounded where $$\Cg(W)=\mat{\int_0^u\cg_1(s)ds\\\int_0^v\cg_2(s)ds}.$$
	
	We simply redefine $x=\myprod{\Cg(W),\Fg}_{L^2(B)}$. So that $\|W\|_{L^1(B)}$'s are uniformly bounded if $\cg_i(x)\ge c|x|$ for all $x\in\RR$ and some $c>0$.

	The uniform bound for $\|W\|_{L^1(B)}$  yields that  of $\|DW\|_{L^2(B)}$ (see \cite[Lemma 5.3.5]{dlebook1}) so that the $BMO$ of weak solutions are uniformly small in small balls if $N=2$. We then have $C^{\ag}$ norm of solutions for some $\ag>0$ are uniformly bounded by the theory in \cite{dlebook1}. This provides the existence of a solution of the system.
	
	However, the above $A$ does not allow us to establish the existence of nontrivial coexistence because $\ma^{21}_v=0$.
	
	We consider first the case
	\beqno{Amat1}A(u,v)=[\ma^{ij}]=\cg(u,v)\mat{a_{11}&a_{12}\\a_{21}&a_{22}}.\eeq

	By contradiction we can generalize \mref{intineq0} and prove several  versions of it. One should also note that \mref{intineq0} is invariant with respect to scaling but the inequalities in the sequel are not. 
	
	Let $\cg(W)$ be a positive continuous function of $W$. Define $\fg(W),\llg(W)$  to be the normalized principal eigenpairs of $-\Div(\cg(W)D\fg)=\llg(W)\fg$. We start with the following simple version.

	\blemm{intineqlem} Suppose that $\fg(W)\in L^{q'}(B)$ for some   $q>1$ if $W\in W^{1,2}(B)$ and $\LLg:L^p(B)\to L^q(B)$ is continuous for some $p<2^*$. Moreover, $\|\fg(W)\|_{L^{q'}(B)}=1$; the unit ball of $L^{q'}(B)$ is compact and $\LLg(f)\ge0$ and $\LLg(f)=0$ iff $f=0$. For any $\eg>0$ there is $C(\eg)$ such that \beqno{intineq}\iidx{B}{|W|^2}\le \eg\iidx{B}{|DW|^2}+C(\eg)\iidx{B}{\LLg(W)\fg(W)}.\eeq  
	\elemm
	
	\bproof By contradiction, we can find a sequence $\{W_n\}$ such that $\|W_n\|_{W^{1,2}}$'s are bounded (so that $W_n\to W\ne0$  strongly in $L^p(B)$) and 
	$$\iidx{B}{\LLg(W_n)\fg(W_n)}\to 0.$$ But this impossible as  $\{\fg(W_n)\}$ converges weakly in $L^{q}(B)$ (or the unit ball of $L^{q'}(B)$ is compact?) to some $\fg>0$ because $\fg(W_n)>0$. \eproof
	
	The inequality \mref{intineq} with $\cg=1$ is the key point in the proof of \cite[Lemma 5.3.5]{dlebook1} when one test the system to obtain the bound for $L^1$ norm of solutions. Therefore, we get the same result for $A$ being given by \mref{Amat1}.

	We consider a more general case
	\beqno{Amat2}A(u,v)=[\ma^{ij}]=\mat{a_{11}(u)\cg_1(u,v)&a_{12}(v)\cg_1(u,v)\\a_{21}(u)\cg_2(u,v)&a_{22}(v)\cg_2(u,v)}.\eeq
	
	For $i=1,2$ define $\fg_i(W),\llg_i(W)$  to be the principal eigenpairs of $-\Div(\cg_i(W)D\fg)=\llg_i(W)\fg$. We test the $i^{th}$ equation with $\fg_i$ and integrate by parts twice to obtain 
	$$\iidx{B}{\hat{A}(W)\Cg(W)}=\iidx{B}{G(W)W\odot\Fg(W)}$$
	where $\odot$ is the Hadamard product ($[a_{ij}]\odot[b_{ij}]=[a_{ij}b_{ij}]$) and
	$$\hat{A}(W)=[\hat{a}_{ij}] \mbox{ with } \hat{a}_{ij}(t)=\int_0^ta_{ij}(s)ds,\; \mbox{ where $t=u$ if $j=1$ and $t=v$ if $j=2$},$$
	\beqno{Cgdef}\Cg(W)=\mat{\llg_1(W)\fg_1(W)\\\llg_2(W)\fg_2(W)},\; \Fg(W)=\mat{\fg_1(W)\\ \fg_2(W)}.\eeq
	
	Clearly $\|\hat{A}(W)\Cg(W)\|_{L^1(B)}$ is bounded under the condition \beqno{Gunif}\left|\iidx{B}{G(W)W\odot\Fg(W)}\right|\le C.\eeq
	We then use \mref{intineq} with $\LLg(W)=\max\{(\hat{a}_{i1}(u)+\hat{a}_{i2}(v))\llg_i(W)\}$  to obtain a uniform bound for $L^1$ norm of solutions.
	
	\brem{LLg} In order to use \mref{intineq} with $\LLg(W)=(\hat{a}_{i1}(u)+\hat{a}_{i2}(v))\llg_i(W)$, $i=1,2$, these functions should be  nonnegative continuous functions (keep in mind that $\llg_i(W)>0$ and we are assuming $u,v\ge0$) in $L^q(B)$ for some appropriate $q>1$ and vanish only at 0. We note that $\llg_i\sim \cg_i$.
	\erem

	\brem{Amatrem} The entries $a_{ij}$ of matrix $A$ in \mref{Amat2} can be of the form $a_{i1}(u,v)=D_u a_i(u,v)$ and $a_{i2}(u,v)=D_v a_i(u,v)$ for some function $a_i$.
	
	\erem

	{\bf $N=2$ and the general $A,G$:} We need to bound $\|DW\|_{L^2}$. Assume that for some $k\ge0$ we have (for simplicity we suppose that $a_{ij}$'s are constants)
	$$\cg(W)\sim|W|^k,\;G(W)\sim|W|^k.$$
	
	Testing the system with $W$, \mref{keyboundest} is now $$c\iidx{B}{(1+|W|^{k+1})|DW|^2}\le C\iidx{B}{(1+|W|^k)|W|^2}.$$
	
	We consider the following main hypothesis of our results.
	\bdes\item[L)] Assume that $\cg(W)\sim |W|^{k}$ for some $k\ge0$ and $\cg(W)>\cg_0>0$. 
	Let $\llg(W), \fg(W)$ be the principal eigenpairs of the eigenvalue problem $-\Div(\cg(W) D\fg(W))=\llg(W) \fg(W)$. For any $p,s$ such that $1\le p<2^*$ and $1<s<\frac{Np(k+2)}{(N-1)p(k+2)-2Nk}$, we normalize
	$\|\fg(W)\|_{L^s(B)}=1$.
	
	Also, $\LLg(W)$ is a function satisfying $0\le \LLg(W)\le C|W|^{L_*}+C$ for any $L_*<\frac{2kN+p(k+2)}{2N}$ and constant $C$. $\LLg(f)=0$ iff $f=0$.

	\edes
	
	We then have  \mref{intineq} in the following general result.
	
	\btheo{intlem1} Suppose L) then for any $\eg>0$ there is $C(\eg)>0$ such that
	
	\beqno{intineq1}\iidx{B}{|W|^{k+2}}\le \eg\iidx{B}{|W|^k|DW|^2}+C(\eg)\iidx{B}{\LLg(W)\fg(W)}\eeq for all $W$ such that $D(W^{\frac{k}{2}+1})\in L^2(B)$.  
	\etheo
	\bproof
	The proof is similar to that of \reflemm{intineqlem}. Let us discuss the assumption L) to see thatwe can carry out the contradiction argument in the proof of that lemma to obtain \mref{intineq1}.

	Consider the equation $-\Div(\cg(W) D\fg)=\llg(W) \fg$. Since $\cg(W)\sim|W|^k$ and $D(|W|^{\frac{k}{2}+1})\in L^2(B)$, we have $|W|^\frac{k+2}{2}\in L^p(B)$ for any $p<2^*$ so that $\cg(W)\in L^\frac{p(k+2)}{2k}(B)$ (not uniform). Therefore the equation $-\Div(\cg(W) D\fg)=\llg(W) \fg$  makes sense in a weak form if $\cg(W) D\fg\in L^1(B)$ so that we must have that $D\fg\in L^r(B)$ with $\frac{1}{r}=1-\frac{2k}{p(k+2)}$. Thus,
	$\fg\in L^{s}(B)$ for $\frac{1}{s}>\frac{1}{r}-\frac{1}{N}$  by Sobolev's imbedding theorem and $s$ is specified in L).
	We see that we can normalize $\fg$ such that $\|\fg\|_{L^s(B)}=1$.  We then have the compactness (of bounded  balls in $L^s(B)$) in the proof of \reflemm{intineqlem}, where $s=q'$.
	
	We must have that $\LLg(W)\fg(W)\in L^1(B)$. We need $\LLg(W)\in L^{s'}(B)$ with $\frac{1}{s'}<1-\frac{1}{r}+\frac{1}{N}$. That is, $s'>\frac{Np(k+2)}{2kN+p(k+2)}=\frac{2N}{2kN+p(k+2)}\frac{p(k+2)}{2}$. Since $W\in L^\frac{p(k+2)}{2}(B)$, we see that if $\LLg(W)\le C|W|^{L_*}+C$ for any $L_*<\frac{2kN+p(k+2)}{2N}$ and constant $C$ then $\LLg(W)\in L^{s'}(B)$.
	\eproof

	In fact, as in \refrem{LLg}, we can argue similarly to see that we can redefine $\LLg(W)\fg(W)$ in \mref{intineq1} to  the following result for $W=(u,v)\in\RR^2$ (for simplicity).
	
	\bcoro{intlem2} For $i=1,2$ define $\fg_i(W),\llg_i(W)$  to be the principal eigenpairs of $-\Div(\cg_i(W)D\fg)=\llg_i(W)\fg$ as in \reftheo{intlem1}. Consider the matrix $$A(u,v)=[\ma^{ij}]=\mat{a_{11}(u)\cg_1(W)&a_{12}(v)\cg_1(W)\\a_{21}(u)\cg_2(W)&a_{22}(v)\cg_2(W)}.$$
	Define $\hat{A}(W)=[\hat{a}_{ij}]$ with  $\hat{a}_{ij}(t)=\int_0^ta_{ij}(s)ds$ where $t=u$ if $j=1$ and $t=v$ if $j=2$; and $$\hat{A}_\LLg(W)=\mat{(\hat{a}_{11}(u)+\hat{a}_{12}(v))\llg_1(W)\\(\hat{a}_{21}(u)+\hat{a}_{22}(v))\llg_2(W)},\; \Fg_z(W)=\mat{z_1\fg_1(W)\\ z_2\fg_2(W)}\quad z_1,z_2>0.$$
	
	Assume that  $0\le\hat{a}_{i1}(u)+\hat{a}_{i2}(v)\sim |W|^a$; $k>\max\{0,\frac{a(N-2)-2}{N-1}\}$.

	Then for all $W=(u,v)$ such that $D(W^{\frac{k}{2}+1})\in L^2(B)$ \beqno{intineq2}\iidx{B}{|W|^{k+2}}\le \eg\iidx{B}{|W|^k|DW|^2}+C(\eg)\iidx{B}{\myprod{\hat{A}_\LLg(W),\Fg_z(W)}}.\eeq 
	\ecoro
	
	The proof is obvious. Replacing $\LLg(W)\fg(W)$ in L) by $\myprod{\hat{A}_\LLg(W),\Fg(W)}$, we see that if  $\hat{a}_{i1}(u)+\hat{a}_{i2}(v)\sim |W|^a$ then $|\hat{A}_\LLg(W)|\sim |W|^{k+a}$. We need $k+a<L_*<k+\frac{k+2}{N-2}$ or $k>\max\{0,\frac{a(N-2)-2}{N-1}\}$.

	{\bf Systems:}
	Up to this point, our results hold for general $W$. Now, we consider solutions $W=(u,v)$ of ($A(W)$ does not have to be normally elliptic)
	\beqno{system}-\Div(A(W)DW)=G(W)W,\; A(W)=\mat{a_{11}(u)\cg_1(W)&a_{12}(v)\cg_1(W)\\a_{21}(u)\cg_2(W)&a_{22}(v)\cg_2(W)}, \; G(W)=[\mg_{ij}(W)]. \eeq
	
	By testing the $i^{th}$ equation with $z_i\fg_i$ and integrating by parts twice, it is easy to see that
	$$\iidx{B}{\myprod{\hat{A}_\LLg(W),\Fg_z(W)}}=\iidx{B}{\myprod{G(W)W,\Fg_z(W)}}.$$ 
	We want $\|\myprod{\hat{A}_\LLg(W),\Fg_z(W)}\|_{L^1(B)}$ is bounded which is true as the condition \mref{Gunif} is implied by
	\beqno{Gunif1}\iidx{B}{\myprod{G(W)W,\Fg_z(W)}}\quad \mbox{is bounded from above}\eeq for some constant vector $z>0$.  Because $\Fg(W)\in L^s(B)$, the easiest way to verify this is to show that $\|\myprod{z,G(W)W}\|_{L^{s'}(B)}\le C$ with $s'> 1$.  In checking this condition, we want $s'> 1$ small (or $s$ large). When $N=2$, this is the case because $p\in(1,\infty)$.
	
	Since $\myprod{z,G(W)W}$ is usually involved with $W$, the above argument requires a uniform estimate of $\|W\|_{L^{s'}(B)}$ which can be considered to be a weak assumption but still not easy to established. We will present some special cases to show that it can be done.

	{\bf Uniform $L^1$ boundedness:} 
	We now assume that  $\cg(W)\ge\cg_0$ for some constant $\cg_0>0$ in the equation $-\Div(\cg(W)D\fg)=\llg(W)\fg$.
	
	{\em For suitable $k\ge0$ we can choose $\LLg(W)=|W|$ in L). We can do so if $1<\frac{2kN+p(k+2)}{2N}<k+\frac{k+2}{N-2}$ or $k>\max\{0,\frac{N-4}{N-1}\}$. We also choose $p$ such that $\frac{p(k+2)}{p(k+2)-2k}\le 2^*$ and normalize $\|\fg\|_{L^1(B)}=1$. 
		
		Indeed, from the equation of $\fg$, we easily see that $\llg(W)\sim\cg(W)$ and by testing  with $\fg$
		$$\cg_0\iidx{B}{|D\fg|^2}\le \iidx{B}{\llg(W)\fg^2}\le \|\llg(W)\|_{\frac{p(k+2)}{2k}}\|\fg\|^2_{\frac{p(k+2)}{p(k+2)-2k}}\le C\|\fg\|^2_{\frac{p(k+2)}{p(k+2)-2k}}$$ for some constant $C$ as $\|\llg(W)\|_{\frac{p(k+2)}{2k}}\sim\|\cg(W)\|_{\frac{p(k+2)}{2k}}$ if $\|D(|W|^{\frac{k}{2}+1})\|_{L^2(B)}\le C$ (see L)). 
		Because $\frac{p(k+2)}{p(k+2)-2k}\le 2^*$ and $\|\fg\|_{L^1(B)}=1$, we can apply \mref{intineq0} to the term $\|\fg\|^2_{\frac{p(k+2)}{p(k+2)-2k}}$ with $\eg$ sufficiently small and obtain that $\|D\fg(W)\|_{L^2(B)}\le C$ and the compactness needed in the proof of \reftheo{intlem1} ($\{\fg_n\}$ has a convergent subsequence in some $L^r(B)$ if $\|\fg_n\|_{L^1(B)}=1$).

	}
	
	\bcoro{L1coro}
	Assume that  $\cg_1(W)=\cg_2(W)$ in \refcoro{intlem2} (so that $\fg_1=\fg_2=\fg$) and that $\hat{a}_{ij}(t)\ge c_1t^\ag$ for some $\ag,c_1>0$. Furthermore, let $k>\max\{0,\frac{N-4}{N-1}\}$ and choose $p$ such that $\frac{p(k+2)}{p(k+2)-2k}\le 2^*$ and normalize $\|\fg\|_{L^1(B)}=1$. 
	
	Then, for all $W$  we have 
	\beqno{L1int}\iidx{B}{|W|^{k+2}}\le \eg\iidx{B}{|W|^k|DW|^2}+C(\eg)\iidx{B}{|W|\fg(W)}.\eeq

	In addition, for $z=(z_1,z_2)$ with $z_i>0$, $c_2$, $\bg\ge0$ such that $1+\bg< k+\ag$ suppose that $$\myprod{z,G(W)W}\le c_2(u^{1+\bg}+v^{1+\bg}+1)\quad \forall u,v\ge0.$$
	
	Then, for all solutions $W=(u,v)\ge0$ to the system $-\Div(A(W)DW)=G(W)W$ we have that $\|(u+v)^{1+\bg}\fg\|_{L^1(B)}$ is {\em uniformly} bounded.
	
	Finally, if $A(W)$ is normally elliptic then
	$\|W\|_{W^{1,2}(B)}$ (and thus, $\|W\|_{C^\nu(B)}$ for some $\nu>0$) is {\em uniformly} bounded.

	\ecoro
	
	\bproof By the argument before the corollary, we can choose $\LLg(W)=(u+v)$ and obtain \mref{L1int} for general $W$. As $z_1,z_2>0$, $u,v\ge 0$, $\hat{a}_{ij}(t)\ge c_1t^\ag$ and $\llg(W)\sim \cg(W)\sim |W|^k$, we can use Jensen's inequality ($\|\fg\|_{L^1(B)}=1$) and have that 
	$$\iidx{B}{\myprod{\hat{A}_\LLg(W),\Fg_z(W)}}\ge C\iidx{B}{(u+v)^{k+\ag}\fg(W)}\ge C\left(\iidx{B}{(u+v)^{1+\bg}\fg(W)}\right)^\frac{k+\ag}{1+\bg}.$$
	The assumption on $G(W)$ also implies $\myprod{z,G(W)W}\le C_1((u+v)^{1+\bg}+1)$ for some  $C_1>0$. From this and the above inequality, for all solutions $W=(u,v)\ge0$ to the system $-\Div(A(W)DW)=G(W)W$ we then find $C_2>0$ such that
	$$C\left(\iidx{B}{(u+v)^{1+\bg}\fg(W)}\right)^\frac{k+\ag}{1+\bg}\le  C_1\left(\iidx{B}{(u+v)^{1+\bg}\fg(W)}\right)+C_2.$$ 
	
	Because $1+\bg<k+\ag$, we see that $\|(u+v)^{1+\bg}\fg\|_{L^1(B)}\le C_3$, a constant. 
	From this and \mref{L1int}, if $A(W)$ is normally elliptic then we obtain as usual that $\|W\|_{W^{1,2}(B)}$ (and therefore $\|W\|_{C^\ag(B)}$) is {\em uniformly} bounded for all solutions $W\ge0$ to the system $-\Div(A(W)DW)=G(W)W$. The proof is complete. \eproof

	We now present several examples for solutions $W$ to the  $2\times2$ system \mref{system}. For simplicity, we will assume that $N=2$ although we can consider $N>2$ by choosing appropriate $k$. We will take  $A(W)=[a_{ij}\cg_i(W)]$ with $\cg_i(W)=\dg_i+\ag_i u+\bg_i v$ for some $\ag_i,\bg_i\ge0$ (so $k=1$; in fact we can choose $\cg_i(u,v)=\dg_i+\ag_i u^k+\bg_i v^k$ with any $k\ge0$) and
	$$G(W)=\mat{b_1+c_1u & k_1v+l_1\\k_2u+l_2&b_2+c_2v}\mbox{ so that }G(W)W=\mat{b_1u+l_1v+c_1u^2+k_1v^2\\ l_2u+b_2v+c_2v^2+k_2u^2}.$$
	
	{\bf Example 1:} 
	For $z=[z_1,z_2]^T>0$, we have \beqno{GWW}\myprod{z,G(W)W}=z_1(b_1u+l_1v)+z_2(b_2v+l_2u)+(z_2k_2+z_1c_1)u^2+(z_1k_1+z_2c_2)v^2.\eeq Thus, if $u,v\in L^{2s'}(B)$ then we have \mref{Gunif1} and therefore \mref{intineq2} with $\|\myprod{\hat{A}_\LLg(W),\Fg_z(W)}\|_{L^1(B)}$ being bounded. Here, $s'$ can be any number greater than 1.
	
	The assumption on the integrability of $u,v$ can be greatly relaxed in the following example.
	
	{\bf Example 2:} As $z_i>0$, if $z_ik_i+z_jc_j\le0$ when $i\ne j$ then  $$\myprod{z,G(W)W}\le z_1(b_1u+l_1v)+z_2(b_2v+l_2u)$$ by \mref{GWW}. Thus, if $u,v\in L^{s'}(B)$ then we have \mref{intineq2} and $\|\myprod{\hat{A}_\LLg(W),\Fg_z(W)}\|_{L^1(B)}$ is bounded.
	
	In fact, for any values of $k_i,c_i$ we can always find $z_1,z_2>0$ such that $z_ik_i+z_jc_j\ge0$ so that for such $z_1,z_2$ \mref{GWW} implies the above and  $\myprod{z,G(W)W}\le c(|u|+|v|)^2$ for some constant $c$.
	Thus, if $u,v\in L^{2s'}(B)$  then we have that $\|\myprod{\hat{A}_\LLg(W),\Fg_z(W)}\|_{L^1(B)}$ is bounded. Note also that if $b_i,c_i,k_i,l_i>0$ then $G(W)$ is cooperative when $u,v\ge0$. 
	
	\brem{onecomp} We can choose $b_1=l_2=0$ (and $u>0$ so that $G(W)$ is still cooperative) then $\myprod{z_0,G(W)W}\le C|v|^2$. This provides an interesting observation: if one component of $W$ say $v\in L^{2s'}(B)$ then $W=(u,v)\in C^\ag(B)$. \erem
	
	In fact, we can prove (even in more general cases if $k>0$) that the estimate is uniform in the next example. We will exploit \refcoro{intlem2} with $\bg=0<\ag$ by assuming $\hat{a}_{ij}(t)\ge Ct^\ag$. This example is inspired by the results concerning  constant matrices $A$ (diagonal), $G$ where we could take $s'=1$. 
	
	{\bf Example 3:} Choosing $\cg_1(W)=\cg_2(W)\ge\cg_0>0$ and $\cg(W)\sim|W|$ ($k=1$) and assume that $\hat{a}_{ij}\ge ct^{\ag}$, we have from \mref{GWW} that $\bg=1$ in \refcoro{L1coro}. Thus $k+\ag>1+\bg$ if $\ag>1$. \refcoro{L1coro} then implies that
	the integral of $(u+v)^2\fg$ over $B$ is uniformly bounded if $W\ge0$.
	
	In fact, for any nonzero values of $k_i,c_i$ if we can find $z_1,z_2>0$ such that $z_ik_i+z_jc_j\le0$ so that for such $z_1,z_2$ \mref{GWW} implies the above and  $\myprod{z,G(W)W}\le c(|u|+|v|)$ for some constant $c$. Then, the integral of $(u+v)\fg$ over $B$ is uniformly bounded if $W\ge0$
	
	Of course, if $A(W)$ is normally elliptic then for all non-negative solutions $W$ we have 
	$\|W\|_{C^\nu(B)}$ for some $\nu>0$ is {\em uniformly} bounded.

	\section{The parabolic case:}\label{uniparasec} \eqnoset We try to extend the previous argument to parabolic systems. The uniform bounds in this case play an essential role in the global existence problems and dynamics of solution flows.
	
	Consider the parabolic system \beqno{parasys}\left\{\barr{ll}W_t=\Div(\ma(W) DW)+\mg(W)& \mbox{in $\Og\times(0,T)$,}\\ \mbox{Homogeneous Dirichlet or Neumann boundary conditions}&\mbox{on $\partial\Og\times(0,T)$,}\\W(x,0)=W_0(x)&\mbox{in $\Og$.}\earr \right.\eeq Here, $\Og$ is a bounded 2d domain with smooth $\partial\Og$, $T>0$. Following \cite{Am2}, $W_0\in W^{1,p}(\Og)$ for some $p>2$. $\ma(W),\mg(W)$ are $m\times m$ matrix and $m\times1$ vector respectively with sufficiently smooth entries in $W$.
	We assume the usual normal ellipticity (see \cite{Am2}): for some constant $c$ \beqno{normalell}\myprod{\ma(W) DW,DW}\ge \llg(W)|DW|^2,\; \forall W\in W^{1,2}(\Og) \quad \llg(W)\ge c>0.\eeq
	
	We first have  the following result.
	\btheo{parNis2} Consider the parabolic system \mref{parasys}. For any given $k\ge 0$
	suppose that $\llg(W)\sim |W|^k$ and \beqno{GWcond}|\mg(W)|\le C\min\{(|W|^{\frac{k}{2}+2}+1),(|W|^{k+1}+1)\}.\eeq Then, $\|DW(\cdot,t)\|_{L^2(\Og)}$ is uniformly bounded in terms of $\|W(\cdot,t)\|_{L^1(\Og)}$ in any finite time interval. In particular, for any given $W_0\in W^{1,p}(\Og)$ ($p>2$) $W$ exists globally and classical.

	\etheo
	
	We should note that we can allow $N=3$ in obtaining that $\|DW(\cdot,t)\|_{L^2(\Og)}$ is uniformly bounded in terms of $\|W(\cdot,t)\|_{L^1(\Og)}$ in any finite time interval (but $W$ may not be classical).
	
	We omit the term $\bg(W)DW$ for simplicity (see \cite{dleJMAA,dlebook1}). We will see that $\|W\|_{L^1(\Og)}\le C$ for all $t$ if $\llg(W)>c$ (or $\cg$ in the previous argument) for some $c>0$ and there are some simple algebra relations between $\ma,\mg$.
	Concerning the notation, $\llg$ here (\cite{dlebook1}) now plays the role of $\cg$ in the previous discussion. Some details of the argument in \cite{dlebook1} are presented for our proof later on.
	
	Using the fact that $\frac{d}{dt}\iidx{\Og}{|\ma(W) DW|^2}=2\iidx{\Og}{\myprod{\ma(W) DW, D(\ma(W) W_t)}}$ (see \cite[p. 87]{dlebook1}), we  test the system with $\ma(W)W_t$ to obtain
	\beqno{keypara}\barrl{\frac{d}{dt}\iidx{\Og}{|\ma(W) DW|^2}+\iidx{\Og}{\llg(W)|W_t|^2} \le \iidx{\Og}{\myprod{\mg(W),\ma(W) W_t}}\le}{6cm}& \eg\iidx{\Og}{\llg(W)|W_t|^2}+C(\eg)\iidx{\Og}{\frac{1}{\llg(W)}|\ma^T\mg|^2}.\earr\eeq
	
	Choosing $\eg$ small, we also want that \beqno{intineqp}\iidx{\Og}{\frac{1}{\llg}|\ma^T\mg|^2}\le C\iidx{\Og}{\llg|DW|^2}\iidx{\Og}{|\ma DW|^2}+C(\|W\|_{L^1(\Og)}).\eeq If so then we get for $y(t)=\iidx{\Og}{|\ma DW|^2}$ the following Gr\"owall's inequality
	$$y'(t)\le q(t)y(t)+C(\|W\|_{L^1(\Og)}), \quad q(t)=C\iidx{\Og}{\llg(W)|DW|^2}.$$
	
	Concerning $q(t)$ we have the following simple lemma.
	\blemm{qlem} Assume either that $|g(W)|\le C(|W|^{k+1}+1)$ or $|g(W)W|\le C|W|^{2}$ ($g$ can also depend on $t$) and that $\llg(W)\ge c(|W|^k+1)$ for some $k\ge0$. Then, for any $t>0$ 
	$$\int_0^t\iidx{\Og}{\llg(W)|DW|^2}dt\le \|W_0\|_{L^2(\Og)}+tC(\|W\|_{L^1(\Og)}).$$
	\elemm
	
	The proof is elementary. Assume first that $|\mg(W)|\le C(|W|^{k+1}+1)$ (the case  $\myprod{g(W),W}\le c|W|^2$ is similar). 
	Testing \mref{parasys} with $W$ and integrate over $Q$ and using \mref{GWcond},  we apply the well known inequality $\|f\|_{L^2(\Og)}^2\le \eg\|Df\|_{L^2(\Og)}^2+C(\eg,\|f\|_{L^r(\Og)})$ with $f=|W|^{\frac{k}{2}+1}$ and small $r$ to see that $$\iidx{\Og}{g(W)W}\le C\iidx{\Og}{(|W|^{k+2}+|W|)}\le \eg\iidx{\Og}{\llg(W)|DW|^2} +C(\|W\|_{L^1(\Og)}).$$
	Choosing $\eg$ small and integrating the following over $(0,t)$ to obtain the lemma
	\beqno{lambdaDW}\frac{d}{dt}\iidx{\Og}{W^2}+\iidx{\Og}{\llg(W)|DW|^2}\le \eg\iidx{\Og}{\llg|DW|^2} +C(\|W\|_{L^1(\Og)}).\eeq
	
	\brem{GWWsmall} In fact, if $|g(W)|\le C_0|W|^{k+1}+C$, $\llg(W)\ge c_0(|W|^k+1)$ for some $k\ge0$ and $C_0$ is small in terms of $c_0$ then we can drop the term $\|W\|_{L^1(\Og)}$ in the lemma. Indeed, if we test the system with $W$ and use the ellipcity assumption together with $C|W|\le \eg|W|^2+C(\eg)$ then
	$$\frac{d}{dt}\iidx{\Og}{W^2}+\iidx{\Og}{\llg(W)|DW|^2}\le \iidx{\Og}{C_0|W|^{k+2}+\eg|W|^2}+C(\eg).$$ If $C_0,\eg$ are sufficiently then we obtain the assertion.
	
	\erem
	
	We now see that $\int_{0}^tq(s)ds$ is bounded in terms of $\|W_0\|_{L^2(\Og)}$, $\|W\|_{L^1(\Og)}$ and is uniformly bounded in any finite interval $(0,t)$.
	This implies (because $y(0)$ and $c_0=C(\|W\|_{L^1(\Og)})$ are bounded) 
	$$y(t)\le e^{\int_{t_0}^tq(s)ds}\left(y(0)+c_0\right)<\infty.$$
	We then have that $\|DW\|_{L^2(\Og)}$ and the BMO norm of $W$ are bounded (since $N=2$). The first assertion of the Theorem is proved. The second assertion comes from the theory in \cite{dlebook1}.
	In fact, by defining $Y(t)=y(t)+C_0$, we just need to have that $y'(t)\le q(t)y(t)+C_0(q(t)+1)$.

	We need the following elementary to have \mref{intineqp}.
	
	\blemm{intineqnk} If $N\le3$ and $k\ge0$ then for any $\eg>0$ there is $c(\eg,k)$ such that
	$$\iidx{\Og}{|W|^{2k+4}}\le \eg\left(\iidx{\Og}{|W|^k|DW|^2}\right)^2+C(\eg,k)\left(\iidx{\Og}{|W|}\right)^{2k+4}.$$

	\elemm

	\bproof If $1\le q<\frac{2N}{N-2}N$ and $\bg>0$ then we have the well known interpolation inequality $$\left(\iidx{\Og}{|f|^q}\right)^\frac{1}{q}\le \eg\left(\iidx{\Og}{|Df|^2}\right)^\frac{1}{2}+C(\eg,\bg)\left(\iidx{\Og}{|f|^\bg}\right)^\frac{1}{\bg}.$$
	
	If $N\le3$ we can take $q=4$. Let $f=|W|^{\frac{k}{2}+1}$ then $|Df|^2\sim |W|^k|DW|^2$. Then
	$$\iidx{\Og}{|f|^4}\le \eg\left(\iidx{\Og}{|Df|^2}\right)^2+C(\eg,\bg)\left(\iidx{\Og}{|f|^\bg}\right)^\frac{4}{\bg}.$$
	
	The proof is complete for suitable choice of $\bg$.
	\eproof 
	
	\brem{fandL}If $f=0$ on $\partial\Og$ we don't have the term $\|f\|_{L^1(\Og)}$. In fact, we don't $\eg$ small here as in the elliptic case but we do (and $\|f\|_{L^1(\Og)}$ too)to prove that $\int_{t_0}^tq(s)ds<\infty$.   
	The assumption \mref{Aspec} and the choice of $\fg$ here is not general as in \refcoro{L1coro} of the elliptic case. One may wish to have a parabolic version of the condition L) but this very difficult because (for the time being) one has only the results on principal eigenpairs of periodic  parabolic eigenvalue problems. Furthermore, the behavior of $\llg(W)$ in this case is limited and the fact that $W$ is periodic is not clear.
	
	\erem

	As $|\ma DW|\sim\llg(W)|DW|$, $\llg(W)\ge c |W|^k$ and $|\ma DW|^2\ge c|W|^k|DW|^2$ (this requires that $\ma$ is regular, that is $\llg(W)\ge c>0$ for some constant $c$), 
	we obtain \mref{intineqp} if $|\ma^T\mg|^2\le C \llg(W) (|W|^{2k+4}+1)$. That is,  $|\ma^T\mg|\le C(|W|^{\frac{3k}{2}+2}+1)$ or $|\mg|\le C(|W|^{\frac{k}{2}+2}+1)$ for some constant $C$. These are assumed in
	\reftheo{parNis2} and the proof of it is complete.
	
	{\bf Uniform estimates on $(0,\infty)$:} 
	The following result is an important improvement to \reftheo{parNis2} where we could only prove that $\|DW\|_{L^2(\Og)}$ is uniformly bounded (under the assumption that $\|W\|_{L^1(\Og)}$ is uniformly bounded) in any {\em finite} time interval to obtain the global existence. Here, $\|DW\|_{L^2(\Og)}$ is uniformly bounded in $(0,\infty)$ so that we can establish a stronger result on the existence of global attractors in $W^{1,p}(\Og)$. This also entails several important consequences in the dynamics of the corresponding solution flow. The general SKT systems on planar domains is a typical example.

	\bcoro{uniDucoro} Assume as in \reftheo{parNis2} and suppose further that $\llg(W)\ge c|W|$. Then, there is a finite number $C(\|W\|_{L^1(\Og)})$ such that
	\beqno{uniDubound}\iidx{\Og}{\llg(W)|DW|^2}\le \max\{\iidx{\Og}{|\ma(W_0)DW_0|^2}, C(\|W\|_{L^1(\Og)})\}\quad \forall t\in(0,\infty).\eeq
	
	Furthermore, if $\|W\|_{L^1(\Og)}$ is uniformly bounded on $(0,\infty)$ then the dynamical system associated to \mref{parasys} has a global attractor in $W^{1,p}(\Og)$, $p>2$.
	\ecoro
	
	Note that $W_0\in W^{1,p}(\Og)$, $p>2$, so that $W_0$ is bounded on the 2d domain $\Og$.
	
	Again, we should note that \mref{uniDubound} still holds if $N=3$. In addition, if $\llg(W)\ge c_0>0$ then \mref{uniDubound} is valid if $\|W\|_{L^1(\Og)}$ is replaced by $\|W\|_{L^2(\Og)}$.

	\bproof Set $y(t)=\iidx{\Og}{|\ma(W)DW|^2}$. We  consider intervals where $y$ is decreasing or increasing. On a deceasing interval $(t_0,t_1)$, of course $\iidx{\Og}{\llg|DW|^2}\le y(t)\le y(t_0)$.
	
	If $y'(t)\ge 0$, we have from \mref{keypara} that $$\iidx{\Og}{\llg(W)|W_t|^2} \le C(\eg)\iidx{\Og}{\frac{1}{\llg(W)}|\ma^T\mg|^2}.$$ 
	
	For any  $\eg_0>0$, by \reflemm{intineqnk} if $\llg(W)\ge c |W|^k$ and $|\ma^T\mg|^2\le C \llg(W) (|W|^{2k+4}+1)$ (that is,  $|\mg|\le C(|W|^{\frac{k}{2}+2}+1)$ for some constant $C$) 
	we have
	$$\left(\iidx{\Og}{\frac{1}{\llg}|\ma^T\mg|^2}\right)^\frac{1}{2}\le \eg_0\iidx{\Og}{\llg|DW|^2}+C(\eg_0,\|W\|_{L^1(\Og)}).$$ Therefore, we have $$\iidx{\Og}{\llg(W)|W_t|^2} \le \eg_0\iidx{\Og}{\llg|DW|^2}+C(\eg_0,\|W\|_{L^1(\Og)}).$$
	
	Testing the system with $W$ (see \mref{lambdaDW}), we have $$\iidx{\Og}{\llg|DW|^2}\le \left|\iidx{\Og}{W_tW}\right|+C(\|W\|_{L^1(\Og)}).$$
	So that with $C_1=\|\llg^{-1}W^2\|_{L^1(\Og)}^\frac{1}{2}$ which is bounded by $\|W\|_{L^1(\Og)}$ as $\llg(W)\ge c|W|$, by H\"older's inequality and the above two estimates we have some constant $C_2$ depending  on $\|W\|_{L^1(\Og)}$ that
	$$\iidx{\Og}{\llg|DW|^2}\le \left(\iidx{\Og}{\llg|W_t|^2}\right)^\frac{1}{2}\|\llg^{-1}W^2\|_{L^1(\Og)}^\frac{1}{2}+C(\|W\|_{L^1(\Og)})\le \eg_0C_1\iidx{\Og}{\llg|DW|^2}+C_2(\eg_0).$$
	
	We can choose $\eg_0$ small to see that $\iidx{\Og}{\llg|DW|^2}\le C(\eg_0)$ on the increasing intervals of $y(t)$.
	
	Altogether, we obtain the uniform bound \mref{uniDubound}. This gives the uniform bound for $C^\nu(\Og)$ norm of solutions for some $\nu>0$ and completes the proof. \eproof
	
	Another proof is to use the uniform Gronwall inequality (see \cite[Lemma 1.1]{Temam}) for $y'\le gy+h$ for some functions $g,h$ to obtain a uniform on $(0,\infty)$. However we need similar estimates and argument to check the conditions that for some $t_0,r>0$ and finite $a_1,a_2,a_3$
	$$\int_t^{t+r}y(s)ds<a_1,\; \int_t^{t+r}g(s)ds<a_2,\; \int_t^{t+r}h(s)ds<a_3.$$
	
	{\bf On the $L^1$ norm of solutions:}
	Although we can assume a weak assumption that $L^1$ norms of solutions are uniformly bounded but proving this is not a simple matter. As in the elliptic case, we can replace $\|f\|_{L^1(\Og)}$ by $\|f\fg\|_{L^1(\Og)}$ where $\fg>0$ is a suitable function. We can prove that $\|f\fg\|_{L^1(\Og)}$ is uniformly in some special parabolic cases.
	
	Indeed, let $\fg>0$ be the principal eigenfunction of $-\Delta\fg=\llg\fg$ on $\Og$ with $\|\fg\|_{L^1(\Og)}=1$. The 
	\beqno{L1intpara}\iidx{\Og}{|W|^{k+2}}\le \eg\iidx{\Og}{|W|^k|DW|^2}+C(\eg)\iidx{\Og}{|W|\fg}\eeq for all $W$ such that $D(|W|^{k/2+1})\in L^2(\Og)$.
	
	Assume that $W=[u_i]_1^m>0$ satisfies in $\Og\times[0,T]$ the system
	$W_t=\Div(A(W)DW)+G(W)$ with $A(W)=[a_{ij}(u_i)]$ and $G(W)=[g_{ij}(W)]$. Assume that \beqno{Aspec} \sum_{i,j}z_i\int_0^{u_i}a_{ij}(s)ds \ge c_1\sum_i u_i^\ag,\; \sum_{i,j}z_ig_{ij}(W)\le c_2(\sum_i u_i+1) \quad \forall u_i>0\eeq and for {\em some} $z_i,c_1,c_2>0$ and $\ag>1$.
	
	Testing the $i^{th}$ equations with $z_i\fg$ and integrating by parts twice and adding, we obtain
	$$\frac{d}{dt}\iidx{\Og}{\sum_i z_iu_i\fg}=-\llg\iidx{\Og}{\fg\sum_{i,j}z_i\int_0^{u_i}a_{ij}(s)ds}+\iidx{\Og}{\fg\sum_{i,j}z_ig_{ij}(W)}.$$ By the assumption \mref{Aspec}, this and Jensen's inequality (as $\ag>1$) imply some $C>0$ such that
	$$\frac{d}{dt}\iidx{\Og}{\sum_i z_iu_i\fg}\le-C\llg\left(\llg\iidx{\Og}{\fg\sum_{i}z_iu_i}\right)^\ag+C\iidx{\Og}{\fg\sum_{i}z_iu_i}+C.$$
	
	Comparing with the solution of $y'=-C\llg y^\ag+Cy+C$ we see that $\|W\fg\|_{L^1(\Og)}$ is {\em uniformly} bounded if $W>0$ (and because $z_i>0$).
	
	We then have \mref{L1intpara} with $k=\ag$.
	
	{\bf On the $L^2(\Og)$ norm of signed $W$:} We see that proving the uniform boundedness of $\|W\|_{L^1(\Og)}$ is not a clear matter. We can do so for $\|W\fg\|_{L^1(\Og)}$ if $W\ge0$ in some special cases. Proving that $W\ge0$ is not easy for systems. If we can replace $\|W\|_{L^1(\Og)}$ by $\|W\|_{L^2(\Og)}$ then we need not assume that $W\ge0$. Importantly, we can still have the general structure of $\ma$.
	
	We will present some examples when we can do so. Assume in these examples that \beqno{smallGW}|g(W)|\le c_*|W|^{k+1}+C.\eeq
	
	{\em Homogenous Dirichlet boundary conditions:} From \mref{lambdaDW}, if \mref{smallGW} holds for some sufficiently small $c_*>0$ and then we also have
	\beqno{lambdaDW2}\frac{d}{dt}\iidx{\Og}{W^2}+\iidx{\Og}{\llg|DW|^2}\le C\iidx{\Og}{|W|}.\eeq
	Here, we use the imbedding inequality $\|f\|_{L^2(\Og)}^2\le C_*\|Df\|_{L^2(\Og)}^2$ for some constant $C_*$ if $f=0$ on $\partial\Og$. Using Young's inequality, we derive from \mref{lambdaDW2} for some $c>0$ that
	$$\frac{d}{dt}\iidx{\Og}{W^2}+c\iidx{\Og}{|W|^2}\le C|\Og|.$$
	Integrating in $t$ we see that $\|W\|_{L^2(\Og)}$ is uniformly bounded.
	
	{\em Homogeneous Neumann boundary conditions:} This case is a bit more complicated. Let us consider the following simple example. We assume that 
	$g_i(W)\le -c_1u_i+c_0$ for some $c_1>0$, $c_0\ge0$ and all $W=[u_i]_1^m, g(W)=[g_i(W)]_1^m$. 
	We have 
	$$\iidx{\Og}{\sum_ig_i(W)u_i}\le c_*\iidx{\Og}{|W-W_\Og|^{2}}+ \sum_iC((u_i)_\Og),\;\mbox{where } f_\Og=|\Og|^{-1}\iidx{\Og}{f}$$ is the average of $f$ over $\Og$. By the Poincar\'e inequality  $\iidx{\Og}{|W-W_\Og|^{2}}\le C_*\iidx{\Og}{|DW|^{2}}$. Testing the $i^{th}$ equation with $u_i$ and addding the results, if $c_*$ is sufficiently small then for some $c>0$ $$\frac{d}{dt}\iidx{\Og}{W^2}+c\iidx{\Og}{|W|^2}\le \sum_iC((u_i)_\Og).$$ Integrating in $t$, we see that $\|W\|_{L^2(\Og)}$ is bounded in terms of $\|W_0\|_{L^2(\Og)}$ and $(u_i)_\Og$'s.
	
	On the other hand, integrating the $i^{th}$ equation, we have $\frac{d}{dt}(u_i)_\Og \le -c_1(u_i)_\Og+c_0|\Og|$. Since $c_1>0$, we see that $(u_i)_\Og$'s are uniformly bounded and so is $\|W\|_{L^2(\Og)}$.
	
	In the above example, $g_i$ has a linear growth and depends only on $u_i$. By the same argument, we can consider (we assume $k=1$ for simplicity)
	$$g_i(W)=(-c_i-a_i|u_i|+\sum_{j\ne i}b_j|u_j|)u_i.$$
	
	Assume that $a_i,c_i>0$ and $a_i$ large. Then, via Young's inequality,
	$\sum_ig_i(W)u_i\le -c|W|^2$  with $c=\min c_i>0$. We then set $Y=\iidx{\Og}{|W|^2}$ and see that $Y'+cY\le 0$. We conclude that $Y$ is uniformly bounded in terms of $\|W_0\|_{L^2(\Og)}$ as before.
	
	Note that we cannot write $g(W)=\mG(W)W$ for some $\mG(W)$  being cooperative for signed $u_i\ne0$. Let us consider
	$$g_i(W)=(-c_i-a_i|u_i|+\sum_{j\ne i}b_j|u_j|)u_i+\sum_{j\ne i}d_j(|u_j|+1)u_j\mbox{ where $a_i,c_i,d_i>0$ and $a_i,c_i$ large}.$$
	In this case,  $\mG(W)$ is cooperative (as $d_i(|u_i|+1)>0$) however it is not clear that $u_i\ge0$ because the system is not diagonal. The previous argument based on the assumption that $u_i$'s are not negative cannot apply but the assertions still hold because we can still prove that $\|W\|_{L^2(\Og)}$ is uniformly bounded.

	\vspace{2cm}
	{\bf Reaction depends on gradient:} For simplicity, we have assumed that there is no gradient in the reaction. In fact, the same idea can allow $g$ to depend on $W$. We consider the following simple example  to illustrate the ideas
	
	$$W_t=\Div(\ma(W)DW)+\frac{W}{|W|}|DW|^q-f(W) \mbox{ with $1<q<2$, $p=\frac{q}{2-q}>1$,}$$ where $f$ is a vector such that $C(q)|W|^p+Const\le \myprod{f(W),W}$ and  $C(q)$ is the number in the Young inequality $a^qb \le \eg a^2+C(q)b^p$.
	
	We will prove that $q(t)\le C(t)$ if $h(|DW|)\sim |DW|^q$ with $1<q<2$ and $\|W\|_{L^2(\Og)}$ is bounded (assume that for some $W_0$ we can prove somehow that $\|W\|_{L^2(\Og)}\le C(t)$ in $(0,\infty)$). Simply test with $W$ globally and use Young's inequality $|DW|^qW \le \eg|DW|^2+C(q)|W|^p$ for some $\eg<1$. The integral of $C(q)|W|^p$ will be canceled by that of $\myprod{f(W),W}$.

	By \cite{GiaS}, we may have to assume that $\ma(W)$ is bounded (i.e. $k=0$), for any $Q_R=B_R\times (t, t-R^2)$ there is $q_0>1$ such that if $q\le q_0$ then \beqno{higherint}\itQ{Q_R}{|DW|^{2q}}\le C\left(\itQ{Q_{2R}}{|DW|^2}\right)^{q}<C(t)\eeq if $q\sim 1$ and locally. To use higher integrability of $DW$ in $x,t$ above,
	we may have use this result first for $h(|DW|)\sim \eg|DW|^2$ in \cite{GiaS}, which needs crucially that $\ma(W)$ is bounded, and then show that there is some $q_0\in(1,2)$ such that \mref{higherint} holds for $q_0$ and then for all $q\le q_0$.

	Testing with $W_t$ as before and using \mref{higherint}, we have to deal with $|f(W)| |W_t|$  (note that as $q\sim1$  $|f(W)|\sim |W|^{p-1}\le C(|W|+1)$ see \reflemm{qlem} with $k=0$), we get that $\|DW\|_{L^2(\Og)}<\infty$ so that $W$ is BMO (so the solution flow is global (there may be not global attractors) in $W^{1,p_*}(\Og)$ and the gradient is smooth in the interior of $\Og\times (0,T_*)$ when $N=2$ for some $p_*>2$).

	{\bf The effect of cross diffusion} We will present some examples to show that there are essential contrasts between nonlinear diagonal systems and non-diagonal ones although they share the same reaction terms.

	Let $\mu>0$ and $W(x,t)=e^{\mu t}U(x)$ then $W$ solves the evolution system $W_t=\mu W=\Delta W+g(W,t)$ if $-\Delta U+\mu U=f(U)U$ and $g(W,t)=f(U)W=f(e^{-\mu t}W)W$.
	
	By the mountain pass lemma (see \cite[Theorem 3.10]{AR}, here $N\le 3$), there is a nonzero solution $U_0$ to the elliptic system $-\Delta U+\mu U=f(U)U$ (with homogeneous Dirichlet or Neumann boundary conditions), where $f(U)=\mbox{diag}[f_1(U),\ldots,f_m(U)]$ with
	$$f_i(U)=(-a_i|u_i|+\sum_{j\ne i}d_j|u_j|)\mbox{ where $a_i>0$ and $|d_i|$'s are sufficiently small}.$$
	
	We have $\myprod{g(W,t),W}=\myprod{f(U)W,W}\le\myprod{f(e^{-\mu t}W)W,W}\le C|W|^2$ for all $t\ge0$.
	
	We see that the solution $W=e^{\mu t}U_0$ of $W_t=\Delta W+g(W,t)$ (a diagonal system) has $\|W\|_{L^2(\Og)}$ and $\|DW\|_{L^2(\Og)}=e^{\mu t}\|DU_0\|_{L^2(\Og)}$ are  bounded in any finite time interval but NOT on $(0,\infty)$. Therefore, the associated dynamical system does not have a global attractor in $L^2(\Og)$. Meanwhile, for any solution of the cross diffusion system $W_t=\Div(\ma(W) DW)+g(W,t)$ (with $\llg(W)\ge |W|$) we see that (see \reflemm{qlem} and \refcoro{uniDucoro}) $\|W\|_{L^2(\Og)}$ and $\|DW\|_{L^2(\Og)}$ are uniformly bounded  on $(0,\infty)$ (in terms of $\|W(\cdot,0)\|_{L^2(\Og)}$).
	
	The global existence problem of of solutions to diagonal parabolic systems $W_t=\Delta W+g(W,t)$ is pretty classical for suitable $g$. On the other hand, the same problem for cross diffusion systems $W_t=\Div(\ma(W)DW)+g(W,t)$, with $\ma(W)$ being a full matrix, is quite complicated (see \cite{Am2} for local results). Recent global results are reported in \cite{dlebook1} and the results here for $N\le3$ are new.
	
	The above counterexample involves non-autonomous system and we want to deal with autonomous systems.
	
	Let $g(U)=\eg_0|U|^\ag U$, $\ag>0$ and $\eg_0$ is a constant. For $\mu=\mu(t)>0$ if 
	$W_t=\Div(|W|^\ag DW) +g(W)$ and $W_t=\mu'W$ then $-\Div(|W|^\ag DW) +\mu' W=g(W)$.
	
	Making a change of var $\mu v=|W|^\ag W$, then $|W|=|\mu v|^\frac{1}{\ag+1}$ and $W=\mu^\frac{1}{\ag+1} |v|^\frac{-\ag}{\ag+1} v$, we have $\mu Dv=(\ag+1)|W|^\ag DW$  and are led to the equation
	$$-\frac{1}{\ag+1}\Delta v+\frac{\mu'}{\mu}\mu^{\frac{1}{\ag+1}}|v|^\frac{-\ag}{\ag+1}v=\frac{g(\mu^\frac{1}{\ag+1}|v|^\frac{-\ag}{\ag+1}v)}{\mu}.$$
	
	As $g(U)=\eg_0|U|^\ag U$, we have $g(\mu^\frac{1}{\ag+1}|v|^\frac{-\ag}{\ag+1}v)=\eg_0 \mu v$. If we choose $\mu'\mu^\frac{-\ag}{\ag+1}=c_1$ for some positive constant $c_1$ then the above equation is the stationary one $-\frac{1}{\ag+1}\Delta v+c_1|v|^\frac{-\ag}{\ag+1}v=\eg_0 v$. This equation has a nonzero solution $v_0$ via the mountain pass theorem again.
	Let $U_0=|v_0|^\frac{-\ag}{\ag+1}v_0$ ($U_0$ is a nonzero solution of $-\Div(|U|^\ag DU) +c_1U=\eg_0|U|^\ag U$).
	
	We consider the equation $\mu'\mu^\frac{-\ag}{\ag+1}=c_1$ or $\frac{1}{\ag+1}\mu'\mu^{\frac{1}{\ag+1}-1}=\frac{1}{\ag+1}c_1$. It is easy to see that this has a solution $\mu(t)=(\frac{1}{\ag+1}c_1t)^{\ag+1}\to\infty$ as $t\to \infty$.
	
	The above calculations show that the two {\em autonomous} systems (taking $\eg_0=-1$)
	$$W_t=\Div(|W|^\ag DW) +g(W)\mbox{ and }W_t=\Div(\ma(W)DW)+g(W) \mbox{ with $g(U)=-|U|^\ag U$}$$
	exhibit the same phenomena as in the previous example. 
	Namely, the first diagonal (degenerate) system has a solution $W=\frac{1}{\ag+1}c_1tU_0$ with $U_0$ is a nonzero solution of
	$-\Div(|U|^\ag DU) +c_1U=\eg_0|U|^\ag U$ (with  $\ag>0$). This solution has $\|W\|_{L^2(\Og)}=\frac{1}{\ag+1}c_1t\|U_0\|_{L^2(\Og)}$ and $\|DW\|_{L^2(\Og)}=\frac{1}{\ag+1}c_1t\|DU_0\|_{L^2(\Og)}$ are  bounded in any finite time interval but NOT on $(0,\infty)$. For $N=2$ and any given $W_0\in L^2(\Og)$, the associated dynamical system may not be   globally defined. We may have to assume that $W_0$ belongs to smoother space so that its solutions are defined globally but even so the dynamical system does not have a global attractor by the given example. 
	Meanwhile, for any solution of the second cross diffusion system (with $\llg(W)\ge |W|^\ag, \ag>0$) we see that  $\|W\|_{L^2(\Og)}$ and $\|DW\|_{L^2(\Og)}$ are uniformly bounded  on $(0,\infty)$ (in terms of $\|W(\cdot,0)\|_{L^2(\Og)}$ as $\myprod{g(W),W}\le -|W|^2$ if $|W|\ge1$).
	Thus, if $N=2$, its associated dynamical system is globally defined has a global attractor.

	Interestingly, the nonlinear  (degenerate) system $W_t=\Div(|W|^\ag DW) +g(W)$ provides a counterexample of {\em smoothing effect} of linear parabolic systems. The solution $W=\frac{1}{\ag+1}c_1tU_0$ has the same regularity property (in $x$) of its initial data $U_0$. In fact, for some $c_0>0$, we can consider the system $W_t=\Div(|W|^\ag DW)-c_0 W+g(W)$ and see that we have still have a solution $W=c(\mu(t))U_0$ with $\mu$ solving $(\mu' +c_0)\mu^{\frac{-\ag}{\ag+1}}=c_1$. We can see that $c(\mu(t))$ is now bounded but the smoothing effect is not available.

	If $N=2$, cross diffusion system $W_t=\Div(\ma(W)DW)+g(W)$ has this smoothing effect. The solution $W$ is classical if its initial data $W_0\in W^{1,p}(\Og)$ for $p>2$. Thus, the smoothing effect applies.
	
	\vspace{2cm} The following example show the importance of the nonlinearity and growth rate of the diffusion matrices (diagonal or not).
	
	Consider the diagonal semi-linear system
	$W_t=\Div(DW) +C_0|W|W$ We know that $W$ blows up in finite time if $W_t=\Div(DW) +C_0|W|^{2\ag}W$ for some $C_0,\ag>0$, see \cite[Theorems I and II, (III.2) and iv')]{levine}, and the initial data $W(x,0)$ is suitable, see \cite[(III.3)]{levine}. 
	One can see that the abstract settings in \cite{levine} apply to the case $\ma$ is a full {\em constant} matrix (to be precise, the condition (iii) in \cite[Theorem I]{levine} should have been that $\int_0^t(x, \dot{A}x) ds\le0$).

	We will provide some minor modifications to the proof of \cite{levine} to see that the growth rate of $\ma(W)$ is necessary. Let $Au=-\Div(a(t)Du)$ where $a$ is a constant matrix (we can allow $a=a(t)$). Consider the following {\em coupled}  system
	$$\left\{\barr{lll}u_t&=&-Au+{\cal F}(u)+\ccF_1(W)\\W_t&=&\Div(\hat{\ma}(u,W)DW)+\ccF_2(u,W)\earr\right.$$
	Here, we can assume that $\hat{\ma}$ satisfies the ellipticity condition with $\hat{\llg}(u,W)\sim |(u,W)|^\ag$ for some $\ag>0$.  However, the matrix $\ma(u,v)=\mbox{diag}[-A,\hat{\ma}]$ does not satisfy a similar condition with $\llg(u,W)\sim |(u,W)|^\ag$ as we considered previously.  
	
	We change the definition of ${\cal G}$ so that (*) in \cite{levine}  becomes 
	${\cal G}(x,W)=\int_0^1\myprod{{\cal F}(\rg x)+\ccF_1(\rg W),x}d\rg,$  and (**) is now $\frac{d}{dt}{\cal G}(v(t),W)=\myprod{{\cal F}(v(t))+\ccF_1(\rg W),v_t(t)}$.
	
	The proof of \cite[Theorem I]{levine} can be repeated verbatim to the subsystem of $u$ and  shows that the component $u(t)$ can blow up in finite time if ${\cal F}(u) \sim u$ and a suitable $\ccF_1(W)$.
	
	We see that linear cross diffusion (or partly linear) alone is not enough to prevent blow up phenomena because it cannot cause sufficient dissipative effect.
	Meanwhile, by \refrem{GWWsmall} with $k=2\ag$ the nonlinear system $W_t=\Div(\ma(W)DW) +C_0|W|W$ with $C_0>0$ sufficiently small, $\llg(W)\sim |W|^\ag$, and $\ma(W)$ is not degenerate has a global attractor.
	
	\section{Blow-up results}\label{blowupsec}\eqnoset
	
	In this section, we present some  results in \cite{levine,schaffer} on bow-up phenomena of  the {\em scalar} equations. We will generalize the structure of these equation and extend the results to cross diffusion systems. One will see also the contrasts between this phenomenon and the existence of global attractors ($N=2$) obtained earlier when certain parameters are modified. One should remark that blow-up examples were reported (see \cite{JS}) for systems on domains when $N\ge3$ only. There is a blow-up example for $N=2, k=0$ (see \cite{mooney}) but the initial data is discontinuous so that it does not contradict our earlier results.
	
	
	\subsection{A blow-up result for scalar equations:}

	The following theorem generalizes the settings in \cite{schaffer} for  the {\em scalar} equation $u_t=\Delta(a(u))+b(u)$ on a bounded domain $\Og\subset\RR^N$ with initial data $u(0)=u_0\in W^{1,2}(\Og)$. We provide a modified proof, essentially follows \cite{schaffer}, with details in order to have similar results for systems.  
	\btheo{schafferthm} Assume that $a'\ge0$. Let $A(u)=\int_0^ua(s)ds+C$ be a positive potential of $a(u)$. Suppose that there is a function $B$ and \mbox{ for some positive $k<\sqrt{2}$}  \beqno{Bacond}B(u)<a(u)b(u),\; B'(u)=2b(u)a'(u),\; a^2(u)\le k^2a'(u)A(u).\eeq
	
	If $\fg(0):=\iidx{\Og}{A(u_0)}>0$ and $\psi(0):=-\iidx{\Og}{|Da(u_0)|^2}+\iidx{\Og}{B(u_0)}>0$ then $u(t)$ exists on $(0,T)$ with $T<[(c-1)\psi(0)]^{-1}\fg(0)$ where $c=2/k^2>1$.
	
	\etheo
	
	\bproof  Using the equation of $u$ and define $$\fg(t):=\iidx{\Og}{A(u(t))dt}\Rightarrow \fg'(t)=\iidx{\Og}{a(u(t))u_t}=\iidx{\Og}{a(u(t))[\Delta a(u(t))+b(u(t))]}.$$

	Dropping $t$ and integrating by parts, as $B(u)<a(u)b(u)$,  we have
	$$\fg'(t)=-\iidx{\Og}{(|Da(u)|^2+a(u)b(u))} >\left[-\iidx{\Og}{|Da(u)|^2}+\iidx{\Og}{B(u)}\right].$$
	
	We then have $\fg'(t)> \psi(t)$ if we define  $$\psi(t):=-\iidx{\Og}{|Da(u)|^2}+\iidx{\Og}{B(u)}.$$
	
	We have $\psi'(t)=-2\iidx{\Og}{\myprod{Da(u),(Da(u)_t}}+\iidx{\Og}{B_t(u)}$ and by integrating by parts $$\barr{lll}\psi'(t)&=&2\iidx{\Og}{\myprod{\Delta (a(u)),(a(u))_t}}+\iidx{\Og}{B_t(u)}\\&=&2\iidx{\Og}{[u_t-b(u)]a'(u)u_t}+\iidx{\Og}{B_t(u)}\\
	&=&2\iidx{\Og}{a'(u)|u_t|^2}+\iidx{\Og}{[B'(u)-2b(u)a'(u)]u_t}=2\iidx{\Og}{a'(u)|u_t|^2}\earr$$ because $B'(u)=2b(u)a'(u)$. So $\psi'\ge0$. As we impose the condition
	$\psi(0)>0$, we see that $$\psi(t)>0 \mbox{ and so that }\fg'(t)>0\mbox{ for all $t>0$}.$$
	
	Now, by \mref{Bacond}  we have that $a^2(u)\le k^2a'(u)A(u)$ for some $k<\sqrt{2}$ so that $$(\fg'(t))^2\le k^2\left(\iidx{\Og}{\sqrt{a'(u)}u_t\sqrt{A(u)}}\right)^2\le k^2\left(\iidx{\Og}{a'(u)|u_t|^2}\right)\left(\iidx{\Og}{A(u)}\right).
	$$ This is $(\fg')^2\le \frac{k^2}{2}\psi'\fg$. Thus, $\psi\fg'\le\frac{k^2 }{2}\psi'\fg$ or $c\frac{\fg'}{\fg}\le \frac{\psi'}{\psi}$ where $c=2/k^2>1$. Integrating this inequality from 0 to $t$ we obtain
	$$\left(\frac{\fg(t)}{\fg(0)}\right)^c\le \frac{\psi(t)}{\psi(0)}.$$ This implies $\fg'(t)\ge \psi(t)\ge  \psi(0)\fg(0)^{-c}\fg(t)^c$. Further, because $c>1$ and  $\fg(t)>0$,  we have $$0<\fg(t)^{1-c}\le \fg(0)^{1-c}-(c-1)\psi(0)\fg(0)^{-c}t.$$ We must have $t<[(c-1)\psi(0)\fg(0)^{-c}]^{-1}\fg(0)^{1-c}=[(c-1)\psi(0)]^{-1}\fg(0)$. The proof is complete. \eproof
	
	The above theorem allows a more general structure than \cite{schaffer}. The condition $a^2(u)\le k^2a'(u)A(u)$ for some $k<\sqrt{2}$ is obvious when $a(u)=|u|^{m-1}u$ for some $m\ge1$, this was considered in \cite{schaffer}.
	The above result and argument can be extended to  coupled systems in the sequel.
	
	\brem{drift} We write $cDu=\Div(cu)-\Div(c)u$ so that we can consider the system $u_t=\Delta(a(u))+cDu+b(u)$ and write it as $u_t=\Div(D(a(u))+cu)+b(u)-\Div(c)u$.
	
	\erem

	\subsection{Blow-up results for systems:}
	
	Let $a(u)=[a_i(u)]_1^m$ and $b(u)=[b_i(u)]_1^m$, $a_i(u)$, $b_i[u]$ are functions in $u=[u_i]_1^m$. Consider the system $(u_i)_t=\Delta(a_i(u))+b_i(u)$ that is for some $u_0\in W^{1,2}(\Og,\RR^m)$ \beqno{schaffer-sys}
	(u_i)_t=\Div(\frac{\partial a_i(u)}{\partial u_j}Du_j)+b_i(u), \quad u(0)=u_0.\eeq

	We make use of an apparatus in \cite[(*) and (**)]{levine} (with $\ccF(u)=a(u)$) and define
	$$\fg(t)=\int_0^1\iidx{\Og}{\myprod{a(\rg u(t)),u(t)}}d\rg \Rightarrow \fg'(t)=\iidx{\Og}{\myprod{a(u(t)),u_t}}.$$

	Along a solution $u(t)$ of \mref{schaffer-sys},  we now have $$ \fg'(t)=\iidx{\Og}{\myprod{a(u(t)),u_t}}=\iidx{\Og}{\myprod{a(u(t)),[\Delta a(u(t))+b(u(t))]}}.$$
	
	Dropping $t$ and integrating by parts, if $B$ is a {\em scalar function}  such that $B(u)<\myprod{a(u),b(u)}$ then
	$$\fg'(t)=-\iidx{\Og}{(|Da(u)|^2+\myprod{a(u),b(u)})} >\left[-\iidx{\Og}{|Da(u)|^2}+\iidx{\Og}{B(u)}\right].$$
	
	We then have $\fg'(t)> \psi(t)$ if we define  $$\psi(t):=-\iidx{\Og}{|Da(u)|^2}+\iidx{\Og}{B(u)}.$$
	
	So that $\psi'(t)=-2\iidx{\Og}{\myprod{Da(u),(Da(u)_t}}+\iidx{\Og}{B_t(u)}$ and by integrating by parts $$\barr{lll}\psi'(t)&=&2\iidx{\Og}{\myprod{\Delta (a(u)),(a(u))_t}}+\iidx{\Og}{B_t(u)}\\&=&2\iidx{\Og}{\myprod{[u_t-b(u)],a_u(u)u_t}}+\iidx{\Og}{B_t(u)}\\
	&=&2\iidx{\Og}{\myprod{u_t,a_u(u)u_t}}+\iidx{\Og}{[\myprod{B^T_u(u),u_t}-2\myprod{b(u),a_u(u)u_t}}=2\iidx{\Og}{\myprod{u_t,a_u(u)u_t}}\earr$$ if $B^T_{u}(u)=2(a_u(u))^Tb(u)$. So $\psi'\ge0$ if $a_u(u)$ is positive definite. If we also impose the condition
	$$\psi(0)>0 \mbox{ then } \psi(t)>0 \mbox{ and so that }\fg'(t)>0\mbox{ for all $t>0$}.$$
	
	If    $\myprod{a(u),u}\ge0,\; \myprod{a_u(u)\zeta,\zeta}\ge0 \mbox{ for all $u,\zeta\in\RR^m$}$ and for some positive  $\kappa<1/2$
	$$ \myprod{a(u(t)),u_t}^2\le 2\kappa\myprod{u_t,a_u(u)u_t}\int_0^1\myprod{a(\rg u(t)),u(t)}d\rg$$
	then we have  $$(\fg'(t))^2 =\left(\iidx{\Og}{\myprod{a(u(t)),u_t}}\right)^2\le 2\kappa\left(\iidx{\Og}{\left(\myprod{u_t,a_u(u)u_t}\right)^\frac12\left(\int_0^1\myprod{a(\rg u(t)),u(t)}d\rg\right)^\frac12}\right)^2$$ and so that (via H\"older's inequality and the definition of $\fg$)
	$$(\fg'(t))^2 \le 2\kappa\left(\iidx{\Og}{\myprod{u_t,a_u(u)u_t}}\right)\left(\iidx{\Og}{\int_0^1\myprod{a(\rg u(t)),u(t)}d\rg}\right)=2\kappa\psi'(t)\fg(t)$$
	so that $(\fg')^2\le 2\kappa\psi'\fg$. Thus, $\psi\fg'\le (\fg')^2\le2\kappa\psi'\fg$ or $c\frac{\fg'}{\fg}\le \frac{\psi'}{\psi}$ where $c=\frac{1}{2\kappa}>1$ (Therefore, we need $\kappa<1/2$!) This is exactly the differential inequality considered at the end in the proof of \reftheo{schafferthm}, because $c>1$, we must have $t<[(c-1)\psi(0)\fg(0)^{-c}]^{-1}\fg(0)^{1-c}$.

	We summarize  the above in the following 
	\btheo{schafferthmsys} Assume that there is a {\em scalar function} $B$ in $u\in\RR^m$ such that  $$B(u)<\myprod{a(u),b(u)},\; B^T_{u}(u)=2(a_u(u))^Tb(u).$$
	Suppose that \beqno{anonnegative} \myprod{a(u),u}\ge0,\; \myprod{a_u(u)\zeta,\zeta}\ge0 \mbox{ for all $u,\zeta\in\RR^m$}\eeq and, importantly, if for some positive $\kappa<1/2$ 
	\beqno{kappacond}\myprod{a(u(t)),u_t}^2\le 2\kappa\myprod{u_t,a_u(u)u_t}\int_0^1\myprod{a(\rg u(t)),u(t)}d\rg.\eeq
	
	If $\fg(0):=\dspl{\int}_0^1\iidx{\Og}{\myprod{a(\rg u_0),u_0}}d\rg>0$ and $\psi(0):=-\iidx{\Og}{|Da(u_0)|^2}+\iidx{\Og}{B(u_0)}>0$ then the solution $u(t)$ of \mref{schaffer-sys} exists on $(0,T)$ with $T<[(c-1)\psi(0)]^{-1}\fg(0)$ where $c=\frac{1}{2\kappa}>1$.
	
	\etheo

	We note that if $a(u)\sim |u|^ku$ and $b(u)=C_0|u|^ku$ with $C_0>0$ small then the initial data $u_0$  may not exist for $\psi(0)>0$ so that the blow-up theorem does not hold and we still have the existence of global attractors when $N=2$. One can see that our $\ma$ is more general than \mref{schaffer-sys}.
	
	The condition \mref{kappacond} for general systems is much more restricted than \mref{Bacond} of \reftheo{schafferthm}. If $a(u)=\mbox{diag}[a_1(u_1),\cdots,a_m(u_m)]$ is diagonal (and $b(u)$ is a vector) then we can apply the proof for scalar equations to each equation of the system (define $\fg(t)=\sum_i \int_0^{u_i}a_i(s)ds+C$) and obtain the following
	\bcoro{schaffercoro} Assume that $a_i'\ge0$. Let $A(u)=[A_i(u_i)]_1^m$, $A_i(u_i)=\int_0^{u_i}a_i(s)ds+C$ be a positive potential of $a(u)$. Suppose that there is a function $B$ and \mbox{ for some positive $k<\sqrt{2}$}  \beqno{Bacondcoro}B(u)<\myprod{a(u),b(u)},\; B'(u)=2\myprod{a'(u),b(u)},\; |a(u)|^2\le k^2\myprod{a'(u),A(u)}.\eeq
	
	If $\fg(0):=\iidx{\Og}{\sum_{i=1}^mA_i((u_0))_i}>0$ and $\psi(0):=-\iidx{\Og}{|Da(u_0)|^2}+\iidx{\Og}{B(u_0)}>0$ then $u(t)$ exists on $(0,T)$ with $T<[(c-1)\psi(0)]^{-1}\fg(0)$ where $c=2/k^2>1$.
	
	\ecoro
	
	The last condition $|a(u)|^2\le k^2\myprod{a'(u),A(u)}$ in \mref{Bacondcoro} is simply  $a_i^2(u_i)\le k^2a_i'(u_i),A_i(u_i)$ for all $i$.
	
	{\bf Convective terms:} The effect of convection in blow-up phenomena is very delicate and there are very few results on this matter. Especially when $a$ is diagonal as in \refcoro{schaffercoro} and nonlinear, we do not know any reported result. We can consider the system $u_t=\Delta(a(u))+\mathbf{c} Du+b(u)$, $ \mathbf{c} $ may depend on $u$. We now have
	$ \fg'(t)=\iidx{\Og}{\myprod{a(u(t)),u_t}}=\iidx{\Og}{\myprod{a(u(t)),[\Delta a(u(t))+ \mathbf{c} Du+b(u(t))]}}$. By Young's inequality,
	if $B$ is a {\em scalar function}, $\myprod{a(u),b(u)}\ge0$ and $\cg>1$ such that we have $B(u)+\frac{| a(u)\mathbf{c} |^2}{4(\cg-1)|a_u|^2}<\myprod{a(u),b(u)}$ we then have $\fg'(t)> \psi(t)$ if we define  $$\psi(t):=- \cg\iidx{\Og}{|Da(u)|^2}+\iidx{\Og}{B(u)} < -\iidx{\Og}{|Da(u)|^2}+\iidx{\Og}{\myprod{a(u),\mathbf{c}Du+ b(u)}}.$$
	
	So that $\psi'(t)=-2 \cg\iidx{\Og}{\myprod{Da(u),(Da(u)_t}}+\iidx{\Og}{B_t(u)}$ and by integrating by parts $$\barr{lll}\psi'(t)&=&2 \cg\iidx{\Og}{\myprod{[u_t- \mathbf{c} Du-b(u)],a_u(u)u_t}}+\iidx{\Og}{B_t(u)}\\
	&=&2 \cg\iidx{\Og}{\myprod{u_t,a_u(u)u_t}}+\iidx{\Og}{[\myprod{B^T_u(u),u_t}-2 \cg\myprod{ \mathbf{c} Du,a_u(u)u_t}-2 \cg\myprod{b(u),a_u(u)u_t}]}\earr$$ if $B^T_{u}(u)=2 \cg(a_u(u))^Tb(u)$ and  $(a_u(u))^T \mathbf{c} =0$ (the second condition in \mref{anonnegative} does not force $ \mathbf{c} =0$) then we can go on similarly  to obtain blow-ups. 
	In particular, we need that $$\psi(0):=- \cg\iidx{\Og}{|Da(u_0)|^2}+\iidx{\Og}{B(u_0)}>0.$$
	
	Of course, things are complicated if $a_u$ is not a constant matrix.
	
	The above argument apply to diagonal system $u_t=\Delta(a(u))+\mathbf{c} Du+b(u)$ with convection and $a(u)=\mbox{diag}[a_1(u_1),\cdots,a_m(u_m)]$.
	\bcoro{schafferconv} Assume the notations  as in \refcoro{schaffercoro}. Suppose that there is a function $B$ and \mbox{ for some positive $\cg>1$ and $k<\sqrt{2}$}  \beqno{Bacondconv}B(u)+\frac{| a(u)\mathbf{c} |^2}{4(\cg-1)|a_u|^2}<\myprod{a(u),b(u)},\; B'(u)=2\cg\myprod{a'(u),b(u)},\; |a(u)|^2\le k^2\myprod{a'(u),A(u)}.\eeq
	
	If $\fg(0)>0$ and $\psi(0):=-\cg\iidx{\Og}{|Da(u_0)|^2}+\iidx{\Og}{B(u_0)}>0$ and $(a_u(u))^T \mathbf{c} =\myprod{a_u(u), \mathbf{c} }=0$ then $u(t)$ exists on $(0,T)$ with $T<[(c-1)\psi(0)]^{-1}\fg(0)$ where $c=2/k^2>1$.
	
	\ecoro
	
	The above corollary shows the preventive effect of convection ($\mathbf{c}\ne0$) in blow-up phenomenon. As $\cg>1$,  it is harder to have blow-up occurs when even $-\iidx{\Og}{|Da(u_0)|^2}+\iidx{\Og}{B(u_0)}>0$. We need stronger initial data $-\cg\iidx{\Og}{|Da(u_0)|^2}+\iidx{\Og}{B(u_0)}>0$.

	\subsection{Examples:} Among the conditions of \reftheo{schafferthmsys}, \mref{kappacond} seems to be the most complicated (it is easy to find $B$). Roughly speaking, this condition and the choice of initial data state that if both the diffusion and reaction are comparatively strong then blow-ups happen.  This should be compared with the counterexamples in \cite{JS} (for $N\ge 3$ only) where the solutions disappear while they are still bounded. We will take a close look at \mref{kappacond} below and present cases when this could be verified.
	
	{\bf Degenerate strongly coupled systems:} If $a_i(u)=|u|^{m_i-1}u_i$ for $m_i\ge1$ then it is easy to check \mref{anonnegative}. Regarding \mref{kappacond}, $(a_i(u))_{u_j}=(m_i-1)|u|^{m_i-2}\frac{u_iu_j}{|u|}+\dg_{ij}|u|^{m_i-1}=|u|^{m_i-2}\frac{(m_i-1)u_iu_j+\dg_{ij}|u|^2}{|u|}$and $a_i(u)a_j(u)=|u|^{m_i+m_j-2}u_iu_j$. 
	To obtain \mref{kappacond}, we need $$|u|^{m_i+m_j-2}u_iu_j\le 2\kappa(a_i(u))_{u_j}\int_0^1\myprod{a(\rg u),u}d\rg=2\kappa\sum_{k=1}^m|u|^{m_i+m_k-1}\frac{(m_i-1)u_iu_j+\dg_{ij}|u|^2}{(m_k+1)|u|}$$ because $ \int_0^1\myprod{a(\rg u),u}d\rg=\sum_k\frac{1}{m_k+1}|u|^{m_k+1}\ge0$. Thus, we want that for all $i,j$ 
	\beqno{kcond1}|u|^{m_j-2}u_iu_j\le 2\kappa\sum_{k=1}^m|u|^{m_k-2}\frac{(m_i-1)u_iu_j+\dg_{ij}|u|^2}{(m_k+1)}.\eeq 
	
	If $m_i=m_j=m_0$ and $u_iu_j\ge0$ (this is a difficult problem as the system is strongly coupled) then \mref{kcond1} is $$u_iu_j\le \kappa\sum_k\frac{(m_0-1)}{(m_0+1)}[u_iu_j+\dg_{ij}|u|^2].$$
	Of course, if $1<\kappa m\frac{m_0-1}{m_0+1}$ then the above holds true. To have $\kappa<1/2$, we must have $m\ge3$ and $m_0$ sufficiently large ($m_0>\frac{m+2}{m-2}$, $m$ is the number of equations). Since this argument is valid if we know that $u\ge0$ so that we can only assert that
	
	{\em {\em Nonnegative} solutions  of {\em degenerate} \mref{schaffer-sys} with \underline{appropriate initial data} will blow up if $m\ge 3$ and $m_0$ is large.}
	
	This says that the cross diffusion system is not an appropriate model in biology/ecology for more than 3 components ($m\ge3$ and $m_0$ large, that is the cross diffusion is too strong) because $u$ usually denotes densities so that $u\ge0$ and should not blow up for any initial data.
	
	For example, in the SKT model we have $m_0=2$ so that $m$ (the number of equations) should not greater than 6. This applies to $N\ge2$.
	
	One should compare this with the previous global existence results when $N=2$ we see that we cannot have nonnegative solutions in this case. Thus, this assertion may be valid only for $N\ge3$ as we discuss the counterexamples in \cite{JS} later.
	
	{\bf A nondegenerate system:} If $a_i(u)=\llg_iu_i+|u|^{m_i-1}u_i$ for some $\llg_i>0$, $m_i>1$ then similar calculations show that in order to have \mref{kappacond} we need $$(\llg_i+|u|^{m_i-1})(\llg_j+|u|^{m_j-1})u_iu_j\le 2\kappa\sum_k(\frac{1}{m_k+1}|u|^{m_k+1}+\llg_k|u_k|^2)(|u|^{m_i-2}\frac{(m_i-1)u_iu_j+\dg_{ij}|u|^2}{|u|}+\dg_{ij}\llg_i)$$When $i\ne j$ and assume that $u_iu_j\ge0$ it is
	\beqno{kcond2}(\llg_i+|u|^{m_i-1})(\llg_j+|u|^{m_j-1})\le 2\kappa(m_i-1)\sum_k(\frac{1}{m_k+1}|u|^{m_k+m_i-2}+\llg_k|u|^{m_i-3}|u_k|^2).\eeq
	
	We see that it is harder to obtain blow-up results as the above obviously fails if at $|u|$ small.
	
	{\bf Diagonal systems:}
	If $a_i(u)=|u_i|^{m_i-1}u_i$ then $\myprod{a(u(t)),u_t}^2=(\sum_i|u_i|^{m_i-1}u_i(u_i)_t)^2$ and
	$$\myprod{u_t,a_u(u)u_t}\int_0^1\myprod{a(\rg u(t)),u(t)}d\rg=\left(\sum_im_i|u_i|^{m_i-1}(u_i)_t^2\right) \left(\sum_i\frac{1}{m_i+1}|u_i|^{m_i+1}\right). $$
	Thus, \mref{kappacond} holds if there is $\kappa<1/2$ such that
	$$\sum_{i,j}|u_i|^{m_i-1}u_i|u_j|^{m_j-1}u_j(u_i)_t(u_j)_t\le 2\kappa\sum_{i,j}\frac{m_i}{m_j+1}|u_j|^{m_j+1}|u_i|^{m_i-1}(u_i)_t^2.$$
	
	We can easily prove that $u_i\ge0$ in this diagonal  case so we need
	$$\sum_{i,j}u_i^{m_i}u_j^{m_j}(u_i)_t(u_j)_t\le 2\kappa\sum_{i,j}\frac{m_i}{m_j+1}u_j^{m_j+1}u_i^{m_i-1}(u_i)_t^2.$$
	
	If $m\ge 2$, one can see that this is impossible when either $u_i, u_j$ is small or large. However, blow-up results hold if  we refer to \refcoro{schaffercoro} for diagonal systems (the assumption \mref{kappacond} of \reftheo{schafferthmsys} for full systems does not apply well in this case).
	
	{\bf A triangular system:} We present a blow up example for the system \beqno{trianglesys}\left\{\barr{l}u_t=\Delta(a^*(u))+ [\Div(c_i(v)Dv)]_1^m+b^*(u,v),\\v_t=\Div(d(v)Dv)+f(u,v).\earr\right.\eeq
	
	Here, $a^*(u)=\mbox{diag}[a_1(u_1),\ldots,a_m(u_m)]$, with $a_i,d>0$, and $b^*(u,v)=[b_1(u,v),\ldots,b_m(u,v)]^T$ (we also think of $a^*$ as a column vector $[a_1(u_1),\ldots,a_m(u_m)]^T$). We denote the column vector of size $(m+1)$ $$C(v,Dv)=[\Div(c_1(v)Dv),\ldots,\Div(c_m(v)Dv),0]^T,$$ where $c_i$'s are functions. Of course, \mref{trianglesys} is a triangular system and it is easy to see that we can impose conditions on $f$ so that {\em $\|v\|_{C^2(\Og)}$ is uniformly bounded in $t$}. In the sequel, we will present conditions on $a^*,d,c_i,b_i$ such that the component $u$ can blow up in finite time.

	Let $U=\mat{u\\v}$, $a(U)=\mat{a^*(u)\\d(v)}$ and  $b(U)=\mat{b^*(u,v)\\f(u,v)}$. Denote $A_i(u_i)=\int_0^{u_i}a_i(s)ds$, $\dg(v)=\int_0^{v}d(s)ds$, and $\fg(t)=\sum_i\iidx{\Og}{A_i(u_i)}+\iidx{\Og}{\dg(v)}$. We  have as usual (using \mref{trianglesys} and our notations)
	$$ \fg'(t)=\iidx{\Og}{\myprod{a(U(t)),U_t}}=\iidx{\Og}{\myprod{a(U(t)),[\Delta a(U(t))+ C(v,Dv)+b(U(t))]}}.$$ 
	
	Integrating by parts and using Young's inequality, it is easy to see that
	if $B$ is a {\em scalar function}, $\myprod{a(U),b(U)}\ge0$ and $\cg>1$ such that we have $B(U)+\frac{|\myprod{a(U),C(v,Dv)} |^2}{4(\cg-1)|a_U|^2}<\myprod{a(U),b(U)}$ (this where we need $Dv,\Delta v$ are uniformly bounded as $a,b,B$ are independent of $Dv$!) we then have $\fg'(t)> \psi(t)$ if $$-\iidx{\Og}{|Da(U)|^2}+\iidx{\Og}{\myprod{a(U),[C(v,Dv)+b(U)]}}>\psi(t):=- \cg\iidx{\Og}{|Da(U)|^2}+\iidx{\Og}{B(U)}.$$
	
	So that $\psi'(t)=-2 \cg\iidx{\Og}{\myprod{Da(U),(Da(U)_t}}+\iidx{\Og}{B_t(U)}$ and by integrating by parts $$\barr{l}\psi'(t)=2 \cg\iidx{\Og}{\myprod{[U_t- C(v, Dv)-b(U)],a_U(U)U_t}}+\iidx{\Og}{B_t(U)}\\
	=2 \cg\iidx{\Og}{\myprod{U_t,a_U(U)U_t}}+\iidx{\Og}{[\myprod{B^T_U(U),U_t}-2 \cg(\myprod{ C(v, Dv),a_U(U)U_t}+\myprod{b(U),a_U(U)U_t})]}\earr$$ if $B^T_{U}(U)=2 \cg(a_U(U))^Tb(U)$ and $C(v,Dv)$ is small (recall that $\|v\|_{C^2(\Og)}$ is uniformly bounded) such that $\psi'\ge0$ via Young's inequality. We then go on similarly as before to obtain blow-up in $u$ if its initial data $U_0$ of $U$ is appropriate. 
	In particular,  (but $\|v\|_{C^2(\Og)}$ is still uniformly bounded as we design only $b$ in the subsystem for $u$), we need that $$\psi(0):=- \cg\iidx{\Og}{|Da(U_0)|^2}+\iidx{\Og}{B(U_0)}>0.$$
	
	Note that, as the subsystem for $u$ is diagonal and $\|v\|_{C^2(\Og)}$ is bounded, we can also choose  $b^*$ such that $u\ge0$. This also true when $N=2$.

	{\bf Nondegerate nondiagonal for $N\ge2$ (the size of global attractors and stabilities of $0$):}
	We revisit and consider the case $a_i(u)=(K+|u|^2)^{\frac{m_i-1}{2}}u_i$ for some large $K>0$. We have $\myprod{a(u(t)),u_t}^2=(\sum_i(K+|u|^2)^{\frac{m_i-1}{2}}u_i(u_i)_t)^2$ and the $ij$ entries of $a_u(u)$ is $$[a_u(u)]_{ij}=(m_i-1)(K+|u|^2)^{\frac{m_i-3}{2}}u_iu_j+(K+|u|^2)^{\frac{m_i-1}{2}}\dg_{ij},$$ so that $$\myprod{u_t,a_u(u)u_t}=\sum_{i,j}[(m_i-1)(K+|u|^2)^{\frac{m_i-3}{2}}u_i(u_i)_tu_j(u_j)_t+(K+|u|^2)^{\frac{m_i-1}{2}}\dg_{ij}(u_i)_t(u_j)_t],$$
	$$\int_0^1\myprod{a(\rg u(t)),u(t)}d\rg=\sum_i\int_0^1\left(K+\rg^2|u|^2\right)^\frac{m_i-1}{2}2\rg|u|^2 \frac{u_i^2}{2|u|^2} d\rg =\sum_i\frac{\left(K+|u|^2\right)^\frac{m_i+1}{2}}{m_i+1} \frac{u_i^2}{|u|^2}. $$

	If $m_i=m_0$ for all $i$ then we need some $\kappa\in(0,1/2)$ such that \mref{kappacond} or
	$$(m_0+1)(\sum_iu_i(u_i)_t)^2\le 2\kappa \sum_{i,j}[(m_0-1)u_i(u_i)_tu_j(u_j)_t+(K+|u|^2)\dg_{ij}(u_i)_t(u_j)_t],$$ which is written as
	$$(m_0(1-2\kappa)+1+2\kappa)(\sum_iu_i(u_i)_t)^2\le 2\kappa (K+|u|^2)\sum_{i}(u_i)_t^2.$$

	If the initial data $u_0$ satisfies the condition of \reftheo{schafferthmsys} and  $|u|^2<  K$ for $K$ large then the above holds and we have that $u$ will blow up. Thus, if $u$ exists globally then there must be some $(x,t)$ such that $|u(x,t)|\ge K$.
	
	In particular, when $N=2$, if a solution $u$ with initial data $u_0\in W^{1,p}(\Og)$, $p>2$, as in \reftheo{schafferthmsys} 
	$$\fg(0):=\iidx{\Og}{\frac{\left(K+|u_0|^2\right)^\frac{m_0+1}{2}}{m_0+1} }, $$ $$\psi(0):=-\iidx{\Og}{|Da((K+|u_0|^2)^{\frac{m_0-1}{2}}u_0)|^2}+\iidx{\Og}{(K+|u_0|^2)^{\frac{m_0-1}{2}}\myprod{u_0,b(u_0)}}>0,$$
	and $\|u\|_{L^1(\Og)}\le C(t)$, $\forall t>0$, then there must be some $(x,t)$ such that $t<t_0:=[(c-1)\psi(0)]^{-1}\fg(0)$ where $c=\frac{1}{2\kappa}>1$ and $|u(x,t)|\ge K$.
	
	The relationship between the cross diffusion and reaction is quite subtle. We see that when the cross diffusion is strong, e.g. $m_0$ is large, if the reaction is sufficiently strong then the ratio $\fg(0)/\psi(0)$  is not clear in order to estimate the time $u$ will be large somewhere.
	
	In fact,  the size of the domains $\Og$ also plays a role. By a simple scaling $x=Ry$ and taking $\Og=B_1$, we see that the condition that $\psi(0):=-\iidx{B_1}{|Da(u_0)|^2}+\iidx{B_1}{B(u_0)}>0$ may hold for some $u_0$ but $-\iidx{B_R}{|Da(u)|^2}+\iidx{B_R}{B(u)}\le0$ for all $u$ if $R$ is sufficiently small, due to Poincar\'e's inequality, so that blow-up does not occur. This is due to the diffusion would be in effect if the size of the domain is small, as we expected.
	
	Again, this blow-up result does not contradict our previous global existence as we required that $C_0$ is small so that we may not find such initial data due to Sobolev's inequalities. This leads to discussions of instabilities of zero and sizes of global attractors.
	
	When $N=2$, even we know that global attractors exist (say, if $C_0$ is small or by other conditions) and they are contained in some ball $B_M$ then we see that $K$ cannot be too large, say $K> M$. Otherwise, one can find a solution such that there is some $(x,t)$ such that $|u(x,t)|\ge K>M$. Thus, we have an estimate for the size of  global attractors, namely $M\le K$.
	
	On the other hand, we know that if $K$ is large then 0 is stable. However, in the presence of reactions, if we can find an initial data again satisfies the condition in \reftheo{schafferthmsys} and it will blow up. This initial data can be very small by scaling so that 0 is unstable. 
	
	In particular, if $a(u)= \llg u$ for any $\llg>0$ then 0 is Liapunov stable (in appropriate norm spaces). However, if we introduced the cross diffusion $a(u)=(\llg+|u|^2)^{\frac{m_0-1}{2}}u$, with $m_0>1$, (and appropriate reaction $b(u)\sim |u|^{k+1}$ with $k>1$) then it is not. Indeed, if for every 
	$\epsilon >0$, there exists a $\delta >0$ such that, if $\|u(0)\|<\delta$, then for every $
	t\geq 0$ we have $\|u(t)\|<\epsilon$. Choose $\eg<\llg$, the above discussion provide a solution $u$ such that 
	$\|u(t)\|<\epsilon$ for some $t>0$ so that 0 is not Liapunov stable. Note that we cannot use the diagonal perturbation $a_i(u)=\llg_iu_i+|u|^{m_i-1}u_i$ in the diffusion to have the similar assertion.

	{\bf Blow up in gradient:} Let us discuss the counterexamples for $m=N\ge3$  in \cite{JS}. Consider the vector function
	\beqno{counterU} u=\frac{x}{\sqrt{\kappa(1-t)+|x|^2}}, \quad t\in(0,1), x\in\RR^N,\eeq
	where $\kappa\in (0, 2\frac{N-1}{N-2})$. Let $\Og$ be the unit ball. The result in \cite[Theorem 3.1]{JS} shows that there is an elliptic matrix $A(u)$ (real analytic in $u$) such that \mref{counterU} is a bounded (signed) solution to $$ u_t=\Div(A(u)Du) \mbox{ in $\Og\times(0,1)$}.$$
	
	We see that $Du(x)=\left( \frac{1}{\sqrt{\kappa(1-t)+|x|^2}}-\frac{|x|^2}{\kappa(1-t)+|x|^2}\right)Id$. So the gradient of $u$ blows up  in the interior of $\Og$ (when $|x|^2\sim \kappa(1-t)$ and $t\to 1^-$) but $u$ stays bounded. \cite[Theorem 4.1]{JS} also provides a signed solution which is blown up in $\|\cdot\|_{L^\infty(\Og)}$ ($A$ depends directly on $x,t$). One should note that $u$ is a {\em signed} solution and satisfies {\em nonhomogenous} boundary conditions on $\partial\Og$ ($N>2$). It would be interesting to find a blow-up example of a nonnegative solution or one with homogeneous boundary conditions on $\partial\Og$. This was done for a triangular system \mref{trianglesys} even $N=2$.

	Regarding the BMO norm of solutions, the system considered in \cite{JS} is $u_t=\Div(A(u)Du)$ with $A(u)=[A^{\ag\bg}_{ij}(u)]$ where
	$$A^{\ag\bg}_{ij}(u)= \thg\dg_{ij}\dg_{\ag\bg} +A_{i\ag}(u)A_{j\bg}(u),$$
	$$A_{i\ag}(u)=\frac{(N-1-\thg-|u|^2[1+\frac{\kappa}{2(N-1)}])\dg_{i\ag}+\frac{\kappa}{2(N-1)}u^iu^\ag}{\sqrt{N(N-1-\thg)-[2(N-1-\thg)+\frac{\kappa}{2}]|u|^2-\thg|u|^4}}.$$
	Here, $\kappa\in (0, 2\frac{N-1}{N-2})$ and $\thg\in(0,N-2-\frac{\kappa}{2(N-1)})$.
	
	We see that the above example satisfies the assumptions of \reftheo{parNis2} and \refrem{GWWsmall} with $k=0, \mg(W)=0$ but we cannot have that $\|DW\|_{L^N(\Og)}$ is uniformly bounded for $m=N\ge3$ (this holds for $N=2$ and $m\ge N$ as we proved in \reftheo{parNis2})!
	
	{\em Therefore, for $N\ge3$ the regularity and global existence theory in \cite{dlebook1} needs more ad hoc structural assumptions and different approaches in order to establish the critical condition that the BMO norms of solutions are small in small balls.}
	
	The proof of the BMO norm of a solution $W$ is small in small balls should not rely on the estimates of $\|DW\|_{L^N(\Og)}$. This can be seen via Moser-Trudinger inequality.
	
	\subsection{Global existence and regularity in $\RR^3$ - Thin domains} 
	
	We consider the following system on a 3d domain $\Og$, 
	\beqno{parasysz}\left\{\barr{ll}W_t=\Div(\ma(W) DW)+\mg(W)& \mbox{in $\Og\times(0,T)$,}\\ \mbox{Homogeneous Dirichlet or Neumann boundary conditions}&\mbox{on $\partial\Og\times(0,T)$,}\\W(x,0)=W_0(x)&\mbox{in $\Og$.}\earr \right.\eeq
	
	We then have the following 3d version of \reftheo{parNis2}.
	\btheo{parNis2z} Consider the parabolic system \mref{parasysz}. For any given $k\ge 0$
	suppose that $\llg(W)\sim |W|^k$ and \mref{GWcond}$$|\mg(W)|\le C\min\{(|W|^{\frac{k}{2}+2}+1),(|W|^{k+1}+1)\}.$$ Then, $\|DW(\cdot,t)\|_{L^2(\Og)}$ is uniformly bounded in terms of $\|W(\cdot,t)\|_{L^1(\Og)}$ in any finite time interval.
	
	Moreover, assume that $\|W(\cdot,t)\|_{L^1(\Og)}$ is uniformly bounded and   for any given $W_0\in W^{1,p}(\Og)$ ($p>3$) and sufficiently small $\eg, R>0$,  $B_R\times(0,R)\subset \Og$ if 
	\beqno{N3conda}\dspl{\int_0^{R}} \iidx{B_R}{|D_{x_3}W(x,x_3,t)|^2}dx_3\le \eg R \eeq 
	for  all $t\in(0,T)$ (the maximal time existence of $W$) then $W$ exists globally and classical.

	\etheo

	\bproof As we noted earlier, \reftheo{parNis2} also yields that $\|DW\|_{L^2(\Og)}$ is uniformly bounded. But $DW\in L^2(\Og)$ implies only that $W(\cdot,x_3,t)$ has small  BMO norms on small {\em disks} in $\RR^2$ but this smallness is not uniform in (a.e.) $x_3$. If  $\dspl{\int_0^\dg}\iidx{\Og\cap\{x_3=s\}}{|DW|^2}ds\le M$ for some constant $M$ then there must be $s\in(0,\dg)$ such that $\iidx{\Og\cap\{x_3=s\}}{|DW|^2}\le M/\dg$. Pick $s\in (0,\dg)$ such that for any ball $B_{3,R}$ in $\RR^3$ with $R$ small, we have that $\|W\|_{BMO(\Og\cap B_{3,R}\cap\{x_3=s\})}\le \eg$.
	
	In the sequel,  we will drop the variable $t$ in our calculations for simplicity of presentation. We also write $z\in \RR^3, x\in\RR^2$.
	
	For $Q=B_R\times(0,R)$, $B_R\subset\RR^2$, we apply the formula $f(x_3)\le\int_s^{x_3}|f'(t_3)|dt_3 +f(s)$ to the function $f(x_3)=\frac{1}{R^2}\iidx{B_R}{|W(x,x_3)-W_{B_R}|}$ (the subscript 3 emphasizes that the derivatives are the partial derivative in the 3rd variable) and appropriate uses of H\"older's inequality to see that

	$$\barrl{\frac{1}{R^3}\dspl{\int_Q}|W-W_Q|dz \le \frac{1}{R}\int_0^R\frac{1}{R^2}\iidx{B_R}{|W-W_{B_R}|}dx_3}{1cm}
	&\le \frac{1}{R}\dspl{\int_0^R}\dspl{\int_s^{x_3}}\left|D_{t_3}\frac{1}{R^2}\iidx{B_R}{|W(x,t_3)-W_{B_R}|}dt_3\right|dx_3+\frac{1}{R^2}\iidx{B_R}{|W(x,s)-W_{B_R}|}
	\\
	&\le \frac{C}{R^3}\dspl{\int_0^R}\dspl{\int_s^{x_3}}\iidx{B_R}{|D_{t_3}W|}dt_3dx_3+\frac{1}{R^2}\iidx{B_R}{|W(x,s)-W_{B_R}|}
	\\
	&\le \frac{C}{R^3}\dspl{\int_0^R}\dspl{\int_s^{x_3}}\left( \iidx{B_R}{|D_{t_3}W|^2}\right)^\frac{1}{2}Rdt_3dx_3+\frac{1}{R^2}\iidx{B_R}{|W(x,s)-W_{B_R}|}
	\\
	&\le \frac{C}{R^2}\dspl{\int_0^R}\left(\dspl{\int_s^{x_3}} \iidx{B_R}{|D_{t_3}W|^2}dt_3\right)^\frac{1}{2}|x_3-s|^\frac{1}{2}dx_3+\frac{1}{R^2}\iidx{B_R}{|W(x,s)-W_{B_R}|}
	\\
	&\le \frac{C}{R^\frac{1}{2}}\left(\dspl{\int_s^{x_3}} \iidx{B_R}{|D_{t_3}W|^2}dt_3\right)^\frac{1}{2}+\frac{1}{R^2}\iidx{B_R}{|W(x,s)-W_{B_R}|}
	\earr$$
	
	Thus, the BMO norm in $\RR^3$ of $W$ (one can replace balls/boxes by cylinders) will be small ($\le (C\sqrt{\eg}+\eg$) if  \mref{N3conda} holds. The proof is complete by invoking the theory in \cite{dlebook1}. \eproof
	
	\brem{GiaSrem} Note that we have $\dspl{\int_0^\dg}\iidx{\Og\times\{s\}}{|DW|^2}ds\le M$ but this does not imply  \mref{N3conda}. This condition relaxes the one in \cite{GiaS} when $N=3$ where it was required that $\dspl{\int_0^{R}} \iidx{B_R}{|DW|^2}dt_3\le \eg R$ for $R$ small (here, $|DW|$ is replaced by $|D_{x_3}W|$ which is of course smaller). 
	It is easy to see that the result can be extended to general $N>3$ if we have $W$ has BMO norms small in small disks of some hyper-plane ($\RR^{N-1}$) and \mref{N3conda} is now
	$$ \dspl{\int_0^{R}} \iidx{B_R}{|D_{x_N}W(x,x_N,t)|^2}dx_N\le \eg R \quad\mbox{for sufficiently small $R$ ($B_R\subset 
		\Og\cap\RR^{N-1}\times\{x_N\}$)}.$$
	We also note that the above can be proved  independently of the boundary conditions as it is local.
	
	\erem
	
	On the other hand, the counterexample in \cite{JS}, presented above with nonhomogeneous boundary conditions, does not give that $Du$ (and $D_{x_3}u$) is bounded. If $D_{x_3}W$ is bounded then
	$R\sup_Q|D_{x_3}W|\le \eg$ when $R$ is small. Of course, this implies \mref{N3conda}.

	We present an example where \reftheo{parNis2z} applies. 
	\bcoro{thindomain} Consider \mref{parasysz}. If $\Og=S\times(0,\dg)$ is a slab in $\RR^3$ ($S$ is a bounded domain in $\RR^2$ and $\dg>0$) and we assume homogenous Neumann boundary condition $D_{x_3}W=0$ on the boundary $\partial\Og_0=S$ where $x_3=0$ of $\Og$ ($W$ satisfies homogeneous Dirichlet or Neumann on $\partial\Og\setminus\partial\Og_0$). If $\dg$ is sufficiently small then \mref{N3conda} holds so that for any given $W_0\in W^{1,p}(\Og)$ ($p>3$) \mref{parasysz} has a unique solution which exists globally and classical if  $\|W(\cdot,t)\|_{L^1(\Og)}$ does not blow up.
	\ecoro

	For simplicity, we assume here homogeneous boundary condition. The result holds for more general  boundary conditions, see \refrem{homothindomain} later. The proof is simple. By Lebesgue's theorem
	$$\liminf_{R\to0}\frac{1}{R}\dspl{\int_0^R}\iidx{\Og\cap\{x_3=s\}}{|D_{x_3}W|^2}ds=\iidx{S}{|D_{x_3}W|^2}=0$$ 
	so that  for any given $\eg>0$ there is $R_S>0$ such that $\dspl{\int_0^{R}} \iidx{B_R}{|D_{x_3}W|^2}dt_3\le \eg R$ if $0<R<R_S$ and $B_R\times(0,R)\subset\Og$. If $\Og$ is thin, $\dg\le R_S$ then \mref{N3conda} holds.

	The above argument also applies to more general thin domains, e.g. $\partial\Og_0$ can be any smooth surface. In fact, if $\Og=S\times(0,\dg)$ and there are $\dg_i$, $i=0,\ldots,n$,  such that $\dg_0=0$, $\dg=\dg_n$, $|\dg_i-\dg_{i-1}|$ is sufficiently small and $W_0\in W^{1,p}(\Og)$ such that $D_{x_3}W_0=0$ on $\{x_3=\dg_i\}$ for all $i<n$ then we can solve the problem on each thin slab $S\times(\dg_{i-1},\dg_i)$ and glue the solutions together to obtain a global solution on $\Og$.

	Combining \refrem{GiaSrem} and \refcoro{Nthindomain}, we obtain the following result for arbitrary $N$.
	
	\bcoro{Nthindomain} Consider \mref{parasysz}. If $N>3$,  $\Og=S\times{\stackrel{N-2 \mbox{ times}}{(0,\dg)\times\ldots\times(0,\dg)}}$ is a domain in $\RR^N$ ($S$ is a bounded domain in $\RR^2$, $\dg>0$)  and we assume homogenous Neumann boundary condition  on $S$ and the boundary parts $x_i=0$ for some $i>2$ of $\Og$ ($W$ satisfies homogeneous Dirichlet or Neumann boundary condition on the rest of $\partial\Og$). There is $R_S>0$ such that if $\dg\le R_S$ then for any given initial data $W_0\in W^{1,p}(\Og)$ ($p>N$) \mref{parasysz} has a unique solution which exists globally and classical if $\|W(\cdot,t)\|_{L^1(\Og)}$ does not blow up.

	\ecoro 
	
	The proof is by induction. We easily see that we need only to prove that the BMO norm in $\RR^N$ of $W$ is small in small balls, this was done in \reftheo{parNis2z} for $N=3$ by using \reftheo{parNis2} (which is true only for $N\le3$). By \refrem{GiaSrem}, if we can prove \mref{N3conda} for $N>3$ then one of the two essential ingredients of the proof of \reftheo{parNis2z} is done. This is established easily as we assume homogeneous Neumann on the part of $\partial\Og$ in the $N$th direction. 
	
	The second ingredient is that the BMO norm of $W$ is small in small balls of some hyper-plane $\RR^{N-1}$ of $\Og$. We are no longer able to use \reftheo{parNis2} (which is true only for $N\le3$) for general $N$. However, $W$ is a solution to a similar  system on any hyper-plane with initial data $W_0$ and boundary conditions. Apply \refcoro{thindomain} in finite number of times we see that $W$ is classical and global so that the BMO norm of $W$ is small in small balls of the hyper-planes of $\Og$. We complete the proof.

	In \refcoro{Nthindomain} the boundary part $\Og=S\times{\stackrel{N-3 \mbox{ times}}{(0,\dg)\times\ldots\times(0,\dg)}}$, $N>4$, is a thin domain itself in $\RR^{N-1}$. We can remove this restriction if we can prove that the solution is classical in any hyper-plane of dimension $N-1$ by other means.

	As in \refrem{GiaSrem}, the  establishment of \mref{N3conda} is independent of the boundary conditions as it is local. We can consider periodic (antiperiodic) Dirich let or Neumann boundary conditions in the $x_3$ ($x_N$) direction. Indeed, suppose that $W$ satisfies periodic (resp., antiperiodic) Neumann condition we have $DW(x,0,t)=DW(x,\dg,t)$ (resp., $DW(x,0,t)=-DW(x,\dg,t)$) for all $x\in S, t>0$ so that
	$$\iidx{\partial\Og\cap\{x_3=0\}}{|D_{x_3}W|^2}=\iidx{\partial\Og\cap\{x_3=\dg\}}{|D_{x_3}W|^2}$$
	so that by Rolle's theorem there must be $s\in[0,\dg]$ such that
	$$\iidx{\Og\cap\{x_3=s\}}{|D_{x_3}W|^2}=0.$$
	
	We now follows the proof of \refcoro{thindomain} to prove \mref{N3conda} if $\dg$ is sufficiently small. We also choose the above value of $s$ to be the one used in the proof of \reftheo{parNis2z}. 
	
	Now, suppose that $W$ satisfies periodic (resp., antiperiodic) Dirichlet condition in $x_3$-direction then we have $W(x,0,t)=W(x,\dg,t)$ (resp., $W(x,0,t)=-W(x,\dg,t)$) for all $x\in S, t>0$ so that we can extend $W$ across the plan $x_3=\dg$ in an odd fashion  (i.e., $W(x,\dg-x_3,t)=-W(x,\dg+x_3,t)$). We then obtain a weak solution which satisfies antiperiodic (resp., periodic) Neumann condition on $S\times(0,2\dg)$ and reduced our problem to the previous one.
	
	The above (odd) reflection argument applies also to general Dirichlet or Neumann boundary conditions on $S$.
	
	To the best of our knowledge, the global existence theory of \mref{parasys} (e.g. \cite{Am2}) is not done when the boundary conditions are periodic (even $N=2$). Meanwhile, the regularity results for {\em bounded weak solutions} in \cite{GiaS} is local and independent of the boundary conditions and \mref{N3conda} is sufficient for those results to hold if $\dg$ is small. If $W$ is a bounded weak solution and $N\le 3$ then $\|W\|_{L^2(\Og)}$ is bounded. This discussion leads to
	the following (the result also holds for $N$ dimension domains as in \refcoro{Nthindomain})
	
	\bcoro{thindomainreg} Consider \mref{parasysz}. If $\Og=S\times(0,\dg)$ is a slab in $\RR^3$ ($S$ is a bounded domain in $\RR^2$ and $\dg>0$) and we assume
	
	\bdes\item[PerBC)] Periodic  (resp., antiperiodic) Dirichlet or Neumann boundary conditions in the $x_3$ directions.
	\edes
	Or more generally,
	\bdes
	\item[NBC)] Dirichlet or Neumann boundary condition  on the boundary $\partial\Og$ in the $x_3$ direction.
	\edes

	If $\dg$ is sufficiently small then a {\em bounded weak} solution $W$ of \mref{parasysz} is H\"older continuous 
	\ecoro
	
	\brem{homothindomain} Of course, the above reflection argument in the proof for PerBC) and a simple change of variables (to have global existence results)  could be applied when $W\ne0$ or $D_{x_3}W\ne0$ on $S$ in the proof of \refcoro{thindomain} so that we can assume non-homogeneous Dirichlet or Neumann boundary conditions and obtain the global existence result for thin domains in $\RR^3$.

	\erem
	
	\brem{blowupthindomain} Comparing with the counterexamples  in $\RR^3$ of \cite{JS} ($m=N=3$) and \refrem{homothindomain}, we see that the thinness of the considered domains in $\RR^3$ is sufficient and necessary. As our global existence results 
	apply to the case $m>2$ is arbitrary, we don't know this thinness condition is necessary.
	
	\erem

	\bibliographystyle{plain}

\end{document}